\documentclass{amsart}
\usepackage[all]{xy}
\usepackage{latexsym}
\usepackage{amssymb}
\usepackage{amsmath}
\usepackage{amsthm}
\usepackage{amscd}
\usepackage{fancyhdr}
\usepackage{fullpage}
\usepackage{mathrsfs}
\input cyracc.def
\xyoption{arc}

\usepackage{color}

\usepackage{tikz}
\usepackage{tikz-cd}
\usetikzlibrary{babel}
\usetikzlibrary{matrix}
\usetikzlibrary{shapes}
\usetikzlibrary{arrows}
\usetikzlibrary{calc,3d}
\usetikzlibrary{decorations,decorations.pathmorphing}
\usetikzlibrary{through}
\tikzset{ext/.style={circle, draw,inner sep=1pt},int/.style={circle,draw,fill,inner sep=1pt},nil/.style={inner sep=1pt}}
\tikzset{exte/.style={circle, draw,inner sep=3pt},inte/.style={circle,draw,fill,inner sep=3pt}}
\tikzset{diagram/.style={matrix of math nodes, row sep=3em, column sep=2.5em, text height=1.5ex, text depth=0.25ex}}
\tikzset{diagram2/.style={matrix of math nodes, row sep=0.5em, column sep=0.5em, text height=1.5ex, text depth=0.25ex}}
\tikzset{every picture/.append style={baseline=-.65ex}}

\def\id{{\mbox{1 \hskip -8pt 1}}}
\newcommand{\sgn}{{\mathit s  \mathit g\mathit  n}}
 \newcommand{\lon}{\longrightarrow}
 \newcommand{\bu}{\bullet}
 
 \newcommand{\rar}{\rightarrow}
 \newcommand{\hook}{\hookrightarrow}

\newcommand{\p}{{\partial}}
\newcommand{\Id}{{\mathrm I\mathrm d}}

\newcommand{\Q}{{\mathbb Q}}
 \newcommand{\Z}{{\mathbb Z}}
 \newcommand{\bS}{{\mathbb S}}

 \newcommand{\R}{{\mathbb R}}
 \newcommand{\N}{{\mathbb N}}
 \newcommand{\K}{{\mathbb K}}

\newcommand{\GC}{\mathsf{GC}}
\newcommand{\fGC}{\mathsf{fGC}}
\newcommand{\wh}{\widehat}
\newcommand{\sDer}{\mathsf{der}}

 \newcommand{\ot}{\otimes}

\newcommand{\sF}{{\mathsf F}}

\newcommand{\sT}{{\mathsf{Poly}}}

\newcommand{\Def}{{\mathsf D\mathsf e\mathsf f }}

 \newcommand{\Beq}{\begin{equation}}
 \newcommand{\Eeq}{\end{equation}}
 \newcommand{\Beqr}{\begin{eqnarray}}
 \newcommand{\Eeqr}{\end{eqnarray}}
 \newcommand{\Beqrn}{\begin{eqnarray*}}
 \newcommand{\Eeqrn}{\end{eqnarray*}}
 \newcommand{\Ba}{\begin{array}}
 \newcommand{\Ea}{\end{array}}
 \newcommand{\Bi}{\begin{itemize}}
 \newcommand{\Ei}{\end{itemize}}
 \newcommand{\Bc}{\begin{center}}
 \newcommand{\Ec}{\end{center}}

\newcommand{\fs}{{\mathfrak s}}

\newcommand{\fw}{{\mathfrak w}}
 \newcommand{\f}{{\mathcal O}}
 \newcommand{\cA}{{\mathcal A}}

 \newcommand{\cE}{{\mathcal E}}
 \newcommand{\cF}{{\mathcal F}}

 \newcommand{\caL}{{\mathcal L}}

 \newcommand{\cP}{{\mathcal P}}
 
 \newcommand{\cR}{{\mathcal R}}
 \newcommand{\cS}{{\mathcal S}}

 \newcommand{\ga}{\gamma}
 
 \newcommand{\Ga}{\Gamma}

 \newcommand{\la}{\lambda}

 \newcommand{\Hom}{{\mathrm H\mathrm o\mathrm m}}
 
 \newcommand{\sip}{\smallskip}
 \newcommand{\bip}{\bigskip}
 \newcommand{\mip}{\vspace{2.5mm}}

  \newcommand{\dfcGC}{\mathsf{dfcGC}}
  \newcommand{\dfGC}{\mathsf{dfGC}}

   \newcommand{\dcGC}{\mathsf{dcGC}}

  \newcommand{\OGC}{\mathsf{OGC}}
  \newcommand{\WGC}{\mathsf{dwGC}}

\newcommand{\Ass}{\mathcal{A} \mathit{ss}}

  \newcommand{\swhHoLB}{\mathcal{H}\mathit{olieb}^{_{\atop ^{\bigstar\circlearrowright}}}}
\newcommand{\swHoLB}{\swhHoLB}
\newcommand{\weHoQPcd}{\mathcal{H}\mathit{oqpois}_{c,d}^{_{\atop ^{\circlearrowright}}}}

\newcommand{\weHoQP}{\mathcal{H}\mathit{oqpois}^{_{\atop ^{\circlearrowright}}}}

\newcommand{\wiweHoQP}{\widehat{\mathcal{H}\mathit{oqpois}}^{\Ba{c}\vspace{-1mm}_{\hspace{-2mm}\circlearrowright} \Ea}}

\newcommand{\wHoQPcd}{\mathcal{H}\mathit{oqpois}_{c,d}^{_{\atop ^{\circlearrowright}}}}
\newcommand{\wHoQP}{\mathcal{H}\mathit{oqpois}^{_{\atop ^{\circlearrowright}}}}

\newcommand{\QP}{\mathcal{Q}\mathit{pois}}
 \newcommand{\HoQP}{\mathcal{H}\mathit{oqpois}}

  \newcommand{\LB}{\mathcal{L}\mathit{ieb}}

\newcommand{\LBcd}{\mathcal{L}\mathit{ieb}_{c,d}}

\newcommand{\HoLBcd}{\mathcal{H}\mathit{olieb}_{c,d}}

\newcommand{\wHoLBcd}{\widehat{\mathcal{H}\mathit{olieb}}_{c,d}}

\newcommand{\HoLB}{\mathcal{H}\mathit{olieb}}

 \newcommand{\grt}{\mathfrak{grt}}
 \newcommand{\Der}{\mathrm{Der}}

\let\leq\leqslant
\let\geq\geqslant

\theoremstyle{plain}
\swapnumbers

\newtheorem{prop-def}[theorem]{Proposition-definition}

\newtheorem{main-theorem}{Main~Theorem}[section]
\newtheorem{section-theorem}{Theorem}[section]
\newtheorem{section-corollary}{Corollary}[section]

\theoremstyle{definition}

\begin{document}

\sloppy

 \newenvironment{proo}{\begin{trivlist} \item{\sc {Proof.}}}
  {\hfill $\square$ \end{trivlist}}

\long\def\symbolfootnote[#1]#2{\begingroup%
\def\thefootnote{\fnsymbol{footnote}}\footnote[#1]{#2}\endgroup}

  \title{On deformation quantization of  quadratic\\ Poisson structures}

\author{Anton~Khoroshkin}
\address{Anton~Khoroshkin: Faculty of Mathematics, National Research University Higher School of Economics,
6 Usacheva street, Moscow, Russian Federation}
\email{akhoroshkin@hse.ru}

\author{Sergei~Merkulov}
\address{Sergei~Merkulov:  Department of Mathematics, Luxembourg University, Maison du Nombre, 6 Avenue de la Fonte,
 L-4364 Esch-sur-Alzette,   Grand Duchy of Luxembourg }
\email{sergei.merkulov@uni.lu}


 \begin{abstract}
 We study the deformation complex of the dg wheeled properad of $\Z$-graded quadratic Poisson structures and prove that it is quasi-isomorphic to the even M.\ Kontsevich graph complex. As a first  application we show that the Grothendieck-Teichm\"uller group acts on the genus completion of that wheeled properad faithfully and essentially transitively. As a second application we classify all the universal quantizations of $\Z$-graded quadratic Poisson structures (together with the underlying {\em homogeneous}\, formality maps). In particular we show that two universal quantizations of Poisson structures are equivalent if the agree on generic quadratic Poisson structures.
\end{abstract}
 \maketitle


{\large
\section{\bf Introduction}
}
\label{sec:introduction}

\sip

\subsection{Deformation quantization}  Since the fundamental paper of Maxim Kontsevich
\cite{Ko2} there was a huge progress in our understanding of universal deformation quantizations of {\em generic}\, Poisson structures which lead, in particular, to the complete classification such quantizations
in terms of Drinfeld associators and, moreover, to the computation of the full cohomology of the deformation complex of any formality map from \cite{Ko2} in terms  of the cohomology of the Kontsevich graph complex $\GC_2$ \cite{Ko, Wi1, Do, AM}.

\sip

All the above mentioned results become, however, void when applied to the special
class of {\em linear}\, Poisson structures on a graded vector space $V$ --- the deformation quantization of such structures is unique up to homotopy (and is given by the universal enveloping algebra), the action of the Grothendieck-Teichm\"uller group on linear Poisson structures is hence trivial etc.

\sip

The next level of complexity comes with {\em quadratic}\, Poisson structures on $V$. What happens with the above mentioned beautiful general statements when applied to a {\em generic quadratic}\, Poisson structure? Do we really need associators to quantize them? How rich is the family of homotopy inequivalent universal quantizations of quadratic Poisson structures?

\sip

In this paper we provide answers to both these questions for {\em $\Z$-graded  quadratic Poisson structures}\, on an arbitrary finite-dimensional $\Z$-graded vector space $V$; if $V$ is concentrated in degree zero, say $V=\R^N$, this notion reduces precisely to the ordinary notion of quadratic Poisson structure on $\R^N$. The main motivation to study {\em $\Z$-graded}\,   quadratic Poisson structures comes from a class of so called {\em homogeneous}\, Kontsevich formality maps,
the ones which respect polynomial degrees of polynomial functions on the underlying $\Z$-graded vector space, and which are discussed in more detail below.

\sip

It is worth pointing out right from the start  an important conceptual difference between deformation quantizations of generic $\Z$-graded Poisson structures and of generic $\Z$-graded {\em quadratic}\, ones.
 In the former case the deformation quantization produced a {\em curved}\, $\cA ss_\infty$ algebra structure on the space of formal smooth functions, the curvature term being unavoidable in general; moreover it is that {\em curvature terms}\, which control, rather surprisingly, the homotopy theory of all such universal deformation quantizations and which are  more or less directly connected to the Kontsevich graph complex $\GC_2$ \cite{Do,AM}. By contrast, a deformation quantizations of generic $\Z$-graded {\em quadratic}\, Poisson structure {\em always}\, produces a {\em flat} (or ordinary) $\cA ss_\infty$ algebra structure, with that important for homotopy classifications {\em the curvature term being equal to zero}.  Hence it is not obvious {\em a priori}\ that both quantization theories are controlled essentially by one and the same Kontsevich graph complex $\GC_2$. Which is the main claim of this paper.

\subsection{Homogeneous formality maps versus quadratic Poisson structures} Let $V$ be a finite-dimensional $\Z$-graded vector space, the graded symmetric tensor algebra
$$
\f_V=\bigoplus_{k\geq 0}  \f_V^k, \ \ \  \f_V^k:=\odot^k V, 
$$
can be understood as the algebra of polynomial functions on the linear space $V^*$; it is well-known that the deformation theory of $\f_V$ as an associative algebra is controlled by the Hochschild dg Lie algebra  $\mathsf{Hoch}(\f_V,\f_V)=\oplus_{n\geq 0}\Hom(\ot^n \f_V, \f_V)[1-n]$ generated by multi-differential operators. A {\em M.\ Kontsevich formality map}\, is a $\caL ie_\infty$ quasi-isomorphism of dg lie algebras,
\Beq\label{1: full formality map}
\cF: \sT(V)\lon \mathsf{Hoch}(\f_V,\f_V),
\Eeq
where $\sT(V)=\oplus_{k,l\geq 0}
 \Hom( \wedge^l V,\odot^k V)[1-l]$ is the Lie algebra of polyvector fields on the affine space $V^*$ \cite{Ko2}. If assume that the underlying space $V^*$ is a linear (rather than affine) space, then both sides of the above map come equipped with an extra grading which takes into account the polynomial degrees of formal functions from $\f_V$. We call a quasi-isomorphism above {\em homogeneous}\, if it respects that the polynomial degrees. More precisely, we consider the Lie subalgebra  of $\mathsf{Poly}(\f_V,\f_V)$
$$
\mathsf{Poly}_{(0)}(\f_V,\f_V)=\bigoplus_{k \geq 1}
 \Hom( \wedge^k V,\odot^k V)[1-l],
$$
 and the Lie subalgebra of $\mathsf{Hoch}(\f_V,\f_V)$,
 $$
 \mathsf{Hoch}_{(0)}(\f_V,\f_V)=\bigoplus_{n\geq 1\atop k_1,\ldots, k_n\geq 1}\Hom(\f_V^{k_1}\ot\ldots \ot \f_V^{k_n}, \f_V^{k_1+...+k_n})[1-n]
$$
which are spanned by polyvector fields and, respectively, polydifferential operators which preserve total polynomial degrees of functions. The Hochschild differential preserves the degree grading and, therefore, the cohomology of the homogeneous Hochschild
complex  $\mathsf{Hoch}_{(0)}(\f_V,\f_V)$ is equal precisely to $\mathsf{Poly}_{(0)}(\f_V,\f_V)$. Any
 $\caL ie_\infty$ quasi-isomorphism
\Beq\label{1: homog formality map}
 \cF_{(0)}: \sT_{(0)}(V) \lon  \mathsf{Hoch}_{(0)}(\f_V,\f_V)
 \Eeq
is called a {\em homogeneous}\, formality map. It is not hard to check that the particular formality map (\ref{1: full formality map}) constructed by M.\ Kontsevich in \cite{Ko2} does have this property when restricted to the subspace $\sT_{(0)}(V)\subset \sT(V)$. So the set of such formality maps is non-empty. We prove in this paper that, up to homotopy equivalence, this set is as large as the set of ``full" formality maps (\ref{1: full formality map}) --- it can be identified with the set of Drinfeld's associators.

\sip

 The Maurer-Cartan elements
of the Lie algebra $\sT(V)$ are called {\em $\Z$-graded Poisson structures}. They have a decomposition,
$$
\pi=\sum_{n,m=0}^\infty\pi_n^m,\ \ \ \  \ \pi_n^m\in \Hom(\wedge^n V, \odot^m V)[2-n]
$$
and can be identified with representation of a certain dg wheeled properad ${\HoLB}_{0,1}^{\atop\bigstar\circlearrowright}$ which has been studied in \cite{AM} where it was proven
that its deformation complex is quasi-isomorphic to the Kontsevich graph complex (implying in particular,
that the group of homotopy automorphisms of ${\HoLB}_{0,1}^{\atop\bigstar\circlearrowright}$ can be identified with the famous and mysterious Grothendieck-Teichm\"uller group $GRT$).

\sip

Similarly, Maurer-Cartan elements
$\nu\in  \sT_{(0)}(V)$ are called {\em $\Z$-graded quadratic Poisson structures}\, on $V^*$; they admit a decomposition
$$
\nu=\sum_{n\geq 1} \nu_n, \ \ \ \ \ \nu^n\in \Hom(\wedge^n V, \odot^n V)[2-n]
$$
If the vector space $V$ is concentrated in degree zero, i.e.\ $V\simeq \K^N$, then only the quadratic term $\nu^2: \Hom(\wedge^2 V,\odot^2 V)$ can be non-zero, i.e.\ in that case we recover the standard notion of quadratic Poisson structure. Thus the theory of homogeneous formality maps is the same as the theory of deformation quantizations of  $\Z$-graded quadratic Poisson structures. We will identify each $\Z$-graded quadratic Poisson structure $\nu$ with a representation in $V$ of a certain dg wheeled properad ${\HoQP}_{0,1}^{\atop\circlearrowright}$ which is introduced and studied in this paper. To understand the homotopy theory of a generic {\em finite-dimensional}\, quadratic Poisson structure
one has to compute the cohomology of the derivation complex $\Der(\weHoQP_{0,1})$ of the genus completion of that wheeled prop (see \S 2 for its definition and \S 3 for the computation). The homotopy theory of possibly {\em infinite-dimensional}\, quadratic Poisson structures is controlled by a different complex --- the derivation complex  $\Der(\HoQP_{0,1})$ of the genus completion of the ordinary (unwheeled) prop $\HoQP_{0,1}\subsetneq \weHoQP_{0,1}$ of $\Z$-graded quadratic Poisson structures.
\sip

The dg prop $\HoQP_{0,1}$ and its wheeled closure $\weHoQP_{0,1}$ come naturally  in a {\em family}\, of  props
$\HoQP_{c,d}\subsetneq \weHoQP_{c,d}$ parameterized by two integers $c,d\in \Z$, the props $\HoQP_{c,d}$ being certain  quotients (see \S 2 for details) of the family of properads $\HoLBcd$ whose homotopy theory has been studied in \cite{MW1}; the case $c=0$, $d=1$ corresponds to quadratic Poisson structures while the case $c=d=1$ corresponds to what one might call {\em quadratic homotopy Lie bialgebras}. At present we are not aware of any applications of the latter case so in applications we are most interested in the case $c=1,d=0$.
However we solve the problem of computing the cohomology of derivation complexes for arbitrary values of the parameters $c$ and $d$.

\sip

 Let $\mathsf{FGC}_d$ the Kontsevich graph
complex $\GC_d$ introduced in \cite{Ko} (see also \cite{Wi1} for a detailed description) spanned by {\em not}\, necessary connected graphs with all vertices at least bivalent. This is one of the most important graph complexes in mathematics admitting applications in algebra, geometry, topology and the theory of moduli spaces of algebaric curves (see e.g.\ \cite{Me6} for a review).
It was proven in \cite{Wi1}  that
$$
H^0(\mathsf{FGC_2}^{\geq 2})=\grt,
$$
where $\grt$  is the Lie algebra of the famous Grothendieck-Teichm\"uller group $GRT$ introduced by V.\ Drinfeld  in \cite{Dr}. The main result of this paper states the following.

\subsubsection{\bf Theorem}
\label{1: Theorem on FGC and Der}  (i)  {\em  There is a canonical morphism of dg Lie algebras,
\Beq\label{1: map from FGC to DerQP}
F^\circlearrowright: \mathsf{FGC}^{\geq 2}_{c+d+1} \lon \Der(\weHoQPcd)
\Eeq
which is a quasi-isomorphism. In particular, there is an isomorphism of Lie algebras
$$
H^0(\Der(\weHoQP_{0,1}))=\grt
$$
that is, the Grothendieck-Teichm\"uller group $GRT$ acts up to homotopy faithfully (and essentially transitively) on the completion of the wheeled properad $\weHoQP$ governing {\em finite}-dimensional $\Z$-graded quadratic Poisson structures.
Moreover $H^i(\weHoQP_{0,1})=0$ for $i\leq -1$.}

\sip

(ii) {\em There is an isomorphism of cohomology groups,
$$
H^\bu(\mathsf{FGC}^{\geq 2}_{c+d+1})= H^\bu(\Der(\HoQP_{c,d}))
$$
implying in particular the isomorphism
$$
H^0(\Der(\HoQP_{0,1}))=\K
$$
which implies, that (in sharp contrast to the situation (i)), the prop $\HoQP_{0,1}$ of possibly {\em infinite}-dimensional $\Z$-graded Poisson structures has no homotopy-non-trivial automorphisms other than the standard rescaling of the generators.}

\bip

The proof of this theorem is quite different from proofs of analogous theorems in \cite{MW3} and \cite{AM} for the props of Lie bialgebras. It goes through a new {\em directed weighted}\, version, $\WGC_d$ of the Kontsevich complex $\GC_d$ (the connected version of $\mathsf{FGC}_d$ considered above) which is spanned by directed connected graphs whose vertices $v$ are assigned weights $\fw_v\in \N$ satisfying certain conditions. There is a natural morphism of complex
$$
F: \GC_d^{\geq 2}\lon \WGC_d
$$
which sends a graph from $\GC_d$  into a sum of graphs of directed weighted graphs by assigning directions to the edges in all possible ways and weights to the vertices in all possible admissible ways. One of the central technical results of this paper is Theorem {\ref{4: Theorem on dcGC and dwGC}} which says essentially that this map is a quasi-isomorphism (up to one rescaling class), and hence leads to the conclusions of the above Theorem.

 \subsection{Classification of universal quantizations of $\Z$-graded quadratic Poisson structures} M.\ Kontsevich {\em universal}\,  formality map (\ref{1: full formality map}) given in \cite{Ko2}
 can be equivalently understood as a morphism of dg props \cite{AM}
 $$
\cF:  {c \cA ss}_\infty \lon  \f(\wh{\HoLB}_{0,1}^{\atop\bigstar\circlearrowright}).
$$
where  ${c \cA ss}_\infty$ is the dg free operad of {\em curved}\, strongly homotopy associative algebras,  and $\wh{\HoLB}_{0,1}^{\atop \bigstar\circlearrowright}$ is the genus completed version of the dg prop of $\Z$-graded (not necessarily quadratic) Poisson structure, and   $\f$ is the {\em polydifferential functor}\, from the category of dg props to the category of dg  operads introduced in \cite{MW2} (see \S 2 and 4 for detailed definitions and reminders).
Reversely, any morphism of dg operads $\cF$ as above satisfying certain non-triviality condition  gives us a universal formality map (\ref{1: full formality map}). The deformation complex of any such a universal formality map has been computed in \cite{AM} where it was identified with the graph complex $\mathsf{FGC}_2^{\geq 2}$; the curvature term on the l.h.s.\ of the map $\cF$ played a central role in that classification result.

\sip

Similarly, any {\em universal}\, homogeneous formality map (\ref{1: homog formality map}) can be equivalently understood as a morphism of dg operads,
 $$
\cF_{(0)}:  {\cA ss}_\infty \lon  \f^{\geq 1}(\wh{\HoQP}_{0,1}^{\atop\circlearrowright}).
$$
where  ${\cA ss}_\infty$ is dg free operad of ordinary ({\em flat}) strongly homotopy associative algebras and $\f^{\geq 1}$ is a truncation of the functor $\f$ which does not allow curvature terms.
The M.\ Kontsevich formality morphism \cite{Ko2} when applied to quadratic Poisson structures gives us a concrete example of a universal homogeneous formality map $\cF_{(0)}$.
We prove in \S 5 that {\em the deformation complex of any universal homogeneous formality map $\cF_{(0)}$ is quasi-isomorphic (up to degree shift) to the deformation complex $\Der(\weHoQP_{0,1})$ studied in Theorem {\ref{1: Theorem on FGC and Der}}}. Hence we conclude that there is a canonical morphism of complexes
$$
\mathsf{FGC}_2 \lon \Def\left({\cA ss}_\infty \stackrel{\cF_{(0)}}{\rar}  \f(\wiweHoQP_{0,1})\right)[1]
$$
 which is a quasi-isomorphism. This result implies that the space of all infinitesimal homotopy deformations of $\cF_{(0)}$ can be identified with the Lie algebra of the Grothendieck-Teichm\"uller group,
 $$
H^1\Def\left({\cA ss}_\infty \stackrel{\cF_{(0)}}{\rar}  \f^{\geq 1}(\wiweHoQP_{0,1})\right)=\grt.
$$
It is not hard to see that every such infinitesimal deformation exponentiates to a genuine deformation of $\cF_{(0)}$ implying the second main result of this paper.

\subsubsection{\bf Theorem}\label{1: Class theorem of homog formalities}
(i) {\em The Grothendieck-Teichm\"uller  group $GRT$  acts freely and
transitively on the set $\cS$ of homotopy classes of universal quantizations of $\Z$-graded quadratic Poisson structures. The set  $\cS$ itself can be identified with the set of Drinfeld associators.

(ii) If two universal quantizations of Poisson structures agree on quadratic Poisson structures, then they are homotopy equivalent.}

\sip


 \sip

\subsection{{Structure of the paper}}
In \S 1 we remind a definition of the prop of homotopy Lie bialgebras, introduce the wheeled prop
$\weHoQP_{0,1}$ controlling finite-dimensional $\Z$-graded quadratic Poisson structures, define their derivation complexes and remind the basic notions and facts of the theory of graph complexes.
In \S 2 we prove the first main result of this paper, Theorem 1.2.1. In \S 3 we explain how
a universal homogeneous formality map can be understood as a morphism of dg operads, compute the cohomology of any such a morphism, and finally prove the second main Theorem 1.3.1 of this paper.

\subsection{Some notation} We work in this paper over a field $\K$ of characteristic zero.
If $V=\oplus_{i\in \Z} V^i$ is a graded vector space over $\K$, then
$V[k]$ stands for the graded vector space with $V[k]^i:=V^{i+k}$ and
and $s^k$ for the associated isomorphism $V\rar V[k]$; for $v\in V^i$ we set $|v|:=i$.
 The set $\{1,2, \ldots, n\}$ is abbreviated to $[n]$;  its group of automorphisms is
denoted by $\bS_n$;
the trivial one-dimensional representation of
 $\bS_n$ is denoted by $\id_n$, while its one dimensional sign representation is
 denoted by $\sgn_n$. The cardinality of a finite set $A$ is denoted by $\# A$.

\mip

{\bf Acknowledgement}. The research of A.Kh. was partially supported by the HSE University Basic Research Program. S.M. was partially supported by the University of Luxembourg RSB internal grant. We are grateful to the HSE and the UL for hospitality during our work on this project. We are very grateful to the referees of our paper for their critical remarks and useful suggestions.

\bip

\bip

{
\Large
\section{\bf Derivation complexes of the wheeled prop of\\ $\Z$-graded quadratic Poisson structures}\label{sec:preliminaries}
}

\sip

\subsection{Reminder on  Lie $(c,d)$-bialgebras}
By definition,  $\LBcd$ is a quadratic properad given as the quotient,
$$
\LB_{c,d}:=\cF ree\langle E\rangle/\langle\cR\rangle,
$$
of the free properad generated by an  $\bS$-bimodule $E=\{E(m,n)\}_{m,n\geq 1}$ with
 all $E(m,n)=0$ except
$$
E(2,1):=\id_1\ot \sgn_2^{c}[c-1]=\mbox{span}\left\langle
\Ba{c}\begin{xy}
 <0mm,-0.55mm>*{};<0mm,-2.5mm>*{}**@{-},
 <0.5mm,0.5mm>*{};<2.2mm,2.2mm>*{}**@{-},
 <-0.48mm,0.48mm>*{};<-2.2mm,2.2mm>*{}**@{-},
 <0mm,0mm>*{\circ};<0mm,0mm>*{}**@{},
 <0.5mm,0.5mm>*{};<2.7mm,2.8mm>*{^{_2}}**@{},
 <-0.48mm,0.48mm>*{};<-2.7mm,2.8mm>*{^{_1}}**@{},
 \end{xy}\Ea
=(-1)^{c}
\Ba{c}\begin{xy}
 <0mm,-0.55mm>*{};<0mm,-2.5mm>*{}**@{-},
 <0.5mm,0.5mm>*{};<2.2mm,2.2mm>*{}**@{-},
 <-0.48mm,0.48mm>*{};<-2.2mm,2.2mm>*{}**@{-},
 <0mm,0mm>*{\circ};<0mm,0mm>*{}**@{},
 <0.5mm,0.5mm>*{};<2.7mm,2.8mm>*{^{_1}}**@{},
 <-0.48mm,0.48mm>*{};<-2.7mm,2.8mm>*{^{_2}}**@{},
 \end{xy}\Ea
   \right\rangle
$$
$$
E(1,2):= \sgn_2^{d}\ot \id_1[d-1]=\mbox{span}\left\langle
\Ba{c}\begin{xy}
 <0mm,0.66mm>*{};<0mm,3mm>*{}**@{-},
 <0.39mm,-0.39mm>*{};<2.2mm,-2.2mm>*{}**@{-},
 <-0.35mm,-0.35mm>*{};<-2.2mm,-2.2mm>*{}**@{-},
 <0mm,0mm>*{\circ};<0mm,0mm>*{}**@{},
   <0.39mm,-0.39mm>*{};<2.9mm,-4mm>*{^{_2}}**@{},
   <-0.35mm,-0.35mm>*{};<-2.8mm,-4mm>*{^{_1}}**@{},
\end{xy}\Ea
=(-1)^{d}
\Ba{c}\begin{xy}
 <0mm,0.66mm>*{};<0mm,3mm>*{}**@{-},
 <0.39mm,-0.39mm>*{};<2.2mm,-2.2mm>*{}**@{-},
 <-0.35mm,-0.35mm>*{};<-2.2mm,-2.2mm>*{}**@{-},
 <0mm,0mm>*{\circ};<0mm,0mm>*{}**@{},
   <0.39mm,-0.39mm>*{};<2.9mm,-4mm>*{^{_1}}**@{},
   <-0.35mm,-0.35mm>*{};<-2.8mm,-4mm>*{^{_2}}**@{},
\end{xy}\Ea
\right\rangle
$$
by the ideal generated by the following elements
\Beq\label{R for LieB}
\cR:\left\{
\Ba{c}
\Ba{c}\resizebox{7mm}{!}{
\begin{xy}
 <0mm,0mm>*{\circ};<0mm,0mm>*{}**@{},
 <0mm,-0.49mm>*{};<0mm,-3.0mm>*{}**@{-},
 <0.49mm,0.49mm>*{};<1.9mm,1.9mm>*{}**@{-},
 <-0.5mm,0.5mm>*{};<-1.9mm,1.9mm>*{}**@{-},
 <-2.3mm,2.3mm>*{\circ};<-2.3mm,2.3mm>*{}**@{},
 <-1.8mm,2.8mm>*{};<0mm,4.9mm>*{}**@{-},
 <-2.8mm,2.9mm>*{};<-4.6mm,4.9mm>*{}**@{-},
   <0.49mm,0.49mm>*{};<2.7mm,2.3mm>*{^3}**@{},
   <-1.8mm,2.8mm>*{};<0.4mm,5.3mm>*{^2}**@{},
   <-2.8mm,2.9mm>*{};<-5.1mm,5.3mm>*{^1}**@{},
 \end{xy}}\Ea
 +
\Ba{c}\resizebox{7mm}{!}{\begin{xy}
 <0mm,0mm>*{\circ};<0mm,0mm>*{}**@{},
 <0mm,-0.49mm>*{};<0mm,-3.0mm>*{}**@{-},
 <0.49mm,0.49mm>*{};<1.9mm,1.9mm>*{}**@{-},
 <-0.5mm,0.5mm>*{};<-1.9mm,1.9mm>*{}**@{-},
 <-2.3mm,2.3mm>*{\circ};<-2.3mm,2.3mm>*{}**@{},
 <-1.8mm,2.8mm>*{};<0mm,4.9mm>*{}**@{-},
 <-2.8mm,2.9mm>*{};<-4.6mm,4.9mm>*{}**@{-},
   <0.49mm,0.49mm>*{};<2.7mm,2.3mm>*{^2}**@{},
   <-1.8mm,2.8mm>*{};<0.4mm,5.3mm>*{^1}**@{},
   <-2.8mm,2.9mm>*{};<-5.1mm,5.3mm>*{^3}**@{},
 \end{xy}}\Ea
 +
\Ba{c}\resizebox{7mm}{!}{\begin{xy}
 <0mm,0mm>*{\circ};<0mm,0mm>*{}**@{},
 <0mm,-0.49mm>*{};<0mm,-3.0mm>*{}**@{-},
 <0.49mm,0.49mm>*{};<1.9mm,1.9mm>*{}**@{-},
 <-0.5mm,0.5mm>*{};<-1.9mm,1.9mm>*{}**@{-},
 <-2.3mm,2.3mm>*{\circ};<-2.3mm,2.3mm>*{}**@{},
 <-1.8mm,2.8mm>*{};<0mm,4.9mm>*{}**@{-},
 <-2.8mm,2.9mm>*{};<-4.6mm,4.9mm>*{}**@{-},
   <0.49mm,0.49mm>*{};<2.7mm,2.3mm>*{^1}**@{},
   <-1.8mm,2.8mm>*{};<0.4mm,5.3mm>*{^3}**@{},
   <-2.8mm,2.9mm>*{};<-5.1mm,5.3mm>*{^2}**@{},
 \end{xy}}\Ea
 \ \ , \ \
\Ba{c}\resizebox{8.4mm}{!}{ \begin{xy}
 <0mm,0mm>*{\circ};<0mm,0mm>*{}**@{},
 <0mm,0.69mm>*{};<0mm,3.0mm>*{}**@{-},
 <0.39mm,-0.39mm>*{};<2.4mm,-2.4mm>*{}**@{-},
 <-0.35mm,-0.35mm>*{};<-1.9mm,-1.9mm>*{}**@{-},
 <-2.4mm,-2.4mm>*{\circ};<-2.4mm,-2.4mm>*{}**@{},
 <-2.0mm,-2.8mm>*{};<0mm,-4.9mm>*{}**@{-},
 <-2.8mm,-2.9mm>*{};<-4.7mm,-4.9mm>*{}**@{-},
    <0.39mm,-0.39mm>*{};<3.3mm,-4.0mm>*{^3}**@{},
    <-2.0mm,-2.8mm>*{};<0.5mm,-6.7mm>*{^2}**@{},
    <-2.8mm,-2.9mm>*{};<-5.2mm,-6.7mm>*{^1}**@{},
 \end{xy}}\Ea
 +
\Ba{c}\resizebox{8.4mm}{!}{ \begin{xy}
 <0mm,0mm>*{\circ};<0mm,0mm>*{}**@{},
 <0mm,0.69mm>*{};<0mm,3.0mm>*{}**@{-},
 <0.39mm,-0.39mm>*{};<2.4mm,-2.4mm>*{}**@{-},
 <-0.35mm,-0.35mm>*{};<-1.9mm,-1.9mm>*{}**@{-},
 <-2.4mm,-2.4mm>*{\circ};<-2.4mm,-2.4mm>*{}**@{},
 <-2.0mm,-2.8mm>*{};<0mm,-4.9mm>*{}**@{-},
 <-2.8mm,-2.9mm>*{};<-4.7mm,-4.9mm>*{}**@{-},
    <0.39mm,-0.39mm>*{};<3.3mm,-4.0mm>*{^2}**@{},
    <-2.0mm,-2.8mm>*{};<0.5mm,-6.7mm>*{^1}**@{},
    <-2.8mm,-2.9mm>*{};<-5.2mm,-6.7mm>*{^3}**@{},
 \end{xy}}\Ea
 +
\Ba{c}\resizebox{8.4mm}{!}{ \begin{xy}
 <0mm,0mm>*{\circ};<0mm,0mm>*{}**@{},
 <0mm,0.69mm>*{};<0mm,3.0mm>*{}**@{-},
 <0.39mm,-0.39mm>*{};<2.4mm,-2.4mm>*{}**@{-},
 <-0.35mm,-0.35mm>*{};<-1.9mm,-1.9mm>*{}**@{-},
 <-2.4mm,-2.4mm>*{\circ};<-2.4mm,-2.4mm>*{}**@{},
 <-2.0mm,-2.8mm>*{};<0mm,-4.9mm>*{}**@{-},
 <-2.8mm,-2.9mm>*{};<-4.7mm,-4.9mm>*{}**@{-},
    <0.39mm,-0.39mm>*{};<3.3mm,-4.0mm>*{^1}**@{},
    <-2.0mm,-2.8mm>*{};<0.5mm,-6.7mm>*{^3}**@{},
    <-2.8mm,-2.9mm>*{};<-5.2mm,-6.7mm>*{^2}**@{},
 \end{xy}}\Ea
 \\
(-1)^{cd+c+d} \Ba{c}\resizebox{5mm}{!}{\begin{xy}
 <0mm,2.47mm>*{};<0mm,0.12mm>*{}**@{-},
 <0.5mm,3.5mm>*{};<2.2mm,5.2mm>*{}**@{-},
 <-0.48mm,3.48mm>*{};<-2.2mm,5.2mm>*{}**@{-},
 <0mm,3mm>*{\circ};<0mm,3mm>*{}**@{},
  <0mm,-0.8mm>*{\circ};<0mm,-0.8mm>*{}**@{},
<-0.39mm,-1.2mm>*{};<-2.2mm,-3.5mm>*{}**@{-},
 <0.39mm,-1.2mm>*{};<2.2mm,-3.5mm>*{}**@{-},
     <0.5mm,3.5mm>*{};<2.8mm,5.7mm>*{^2}**@{},
     <-0.48mm,3.48mm>*{};<-2.8mm,5.7mm>*{^1}**@{},
   <0mm,-0.8mm>*{};<-2.7mm,-5.2mm>*{^1}**@{},
   <0mm,-0.8mm>*{};<2.7mm,-5.2mm>*{^2}**@{},
\end{xy}}\Ea
  -
\Ba{c}\resizebox{7mm}{!}{\begin{xy}
 <0mm,-1.3mm>*{};<0mm,-3.5mm>*{}**@{-},
 <0.38mm,-0.2mm>*{};<2.0mm,2.0mm>*{}**@{-},
 <-0.38mm,-0.2mm>*{};<-2.2mm,2.2mm>*{}**@{-},
<0mm,-0.8mm>*{\circ};<0mm,0.8mm>*{}**@{},
 <2.4mm,2.4mm>*{\circ};<2.4mm,2.4mm>*{}**@{},
 <2.77mm,2.0mm>*{};<4.4mm,-0.8mm>*{}**@{-},
 <2.4mm,3mm>*{};<2.4mm,5.2mm>*{}**@{-},
     <0mm,-1.3mm>*{};<0mm,-5.3mm>*{^1}**@{},
     <2.5mm,2.3mm>*{};<5.1mm,-2.6mm>*{^2}**@{},
    <2.4mm,2.5mm>*{};<2.4mm,5.7mm>*{^2}**@{},
    <-0.38mm,-0.2mm>*{};<-2.8mm,2.5mm>*{^1}**@{},
    \end{xy}}\Ea
 + (-1)^{d}
\Ba{c}\resizebox{7mm}{!}{\begin{xy}
 <0mm,-1.3mm>*{};<0mm,-3.5mm>*{}**@{-},
 <0.38mm,-0.2mm>*{};<2.0mm,2.0mm>*{}**@{-},
 <-0.38mm,-0.2mm>*{};<-2.2mm,2.2mm>*{}**@{-},
<0mm,-0.8mm>*{\circ};<0mm,0.8mm>*{}**@{},
 <2.4mm,2.4mm>*{\circ};<2.4mm,2.4mm>*{}**@{},
 <2.77mm,2.0mm>*{};<4.4mm,-0.8mm>*{}**@{-},
 <2.4mm,3mm>*{};<2.4mm,5.2mm>*{}**@{-},
     <0mm,-1.3mm>*{};<0mm,-5.3mm>*{^2}**@{},
     <2.5mm,2.3mm>*{};<5.1mm,-2.6mm>*{^1}**@{},
    <2.4mm,2.5mm>*{};<2.4mm,5.7mm>*{^2}**@{},
    <-0.38mm,-0.2mm>*{};<-2.8mm,2.5mm>*{^1}**@{},
    \end{xy}}\Ea
  + (-1)^{d+c}
\Ba{c}\resizebox{7mm}{!}{\begin{xy}
 <0mm,-1.3mm>*{};<0mm,-3.5mm>*{}**@{-},
 <0.38mm,-0.2mm>*{};<2.0mm,2.0mm>*{}**@{-},
 <-0.38mm,-0.2mm>*{};<-2.2mm,2.2mm>*{}**@{-},
<0mm,-0.8mm>*{\circ};<0mm,0.8mm>*{}**@{},
 <2.4mm,2.4mm>*{\circ};<2.4mm,2.4mm>*{}**@{},
 <2.77mm,2.0mm>*{};<4.4mm,-0.8mm>*{}**@{-},
 <2.4mm,3mm>*{};<2.4mm,5.2mm>*{}**@{-},
     <0mm,-1.3mm>*{};<0mm,-5.3mm>*{^2}**@{},
     <2.5mm,2.3mm>*{};<5.1mm,-2.6mm>*{^1}**@{},
    <2.4mm,2.5mm>*{};<2.4mm,5.7mm>*{^1}**@{},
    <-0.38mm,-0.2mm>*{};<-2.8mm,2.5mm>*{^2}**@{},
    \end{xy}}\Ea
 + (-1)^{c}
\Ba{c}\resizebox{7mm}{!}{\begin{xy}
 <0mm,-1.3mm>*{};<0mm,-3.5mm>*{}**@{-},
 <0.38mm,-0.2mm>*{};<2.0mm,2.0mm>*{}**@{-},
 <-0.38mm,-0.2mm>*{};<-2.2mm,2.2mm>*{}**@{-},
<0mm,-0.8mm>*{\circ};<0mm,0.8mm>*{}**@{},
 <2.4mm,2.4mm>*{\circ};<2.4mm,2.4mm>*{}**@{},
 <2.77mm,2.0mm>*{};<4.4mm,-0.8mm>*{}**@{-},
 <2.4mm,3mm>*{};<2.4mm,5.2mm>*{}**@{-},
     <0mm,-1.3mm>*{};<0mm,-5.3mm>*{^1}**@{},
     <2.5mm,2.3mm>*{};<5.1mm,-2.6mm>*{^2}**@{},
    <2.4mm,2.5mm>*{};<2.4mm,5.7mm>*{^1}**@{},
    <-0.38mm,-0.2mm>*{};<-2.8mm,2.5mm>*{^2}**@{},
    \end{xy}}\Ea
    \Ea
\right.
\Eeq
Note that, when representing elements of all operads and props discussed in this paper as graphs, we tacitly assume that all edges and legs are {\em directed}\, along the flow going from the bottom of the graph to the top.
The minimal resolution $\HoLBcd$ of the properad $\LBcd$ was constructed in \cite{Ko4, MaVo, V}; it
is the free properad generated by the following (skew)symmetric $(m,n)$-corollas of degree $1 +c(1-m)+d(1-n)$
\Beq\label{2: symmetries of HoLiebcd corollas}
\Ba{c}\resizebox{17mm}{!}{\begin{xy}
 <0mm,0mm>*{\circ};<0mm,0mm>*{}**@{},
 <-0.6mm,0.44mm>*{};<-8mm,5mm>*{}**@{-},
 <-0.4mm,0.7mm>*{};<-4.5mm,5mm>*{}**@{-},
 <0mm,0mm>*{};<1mm,5mm>*{\ldots}**@{},
 <0.4mm,0.7mm>*{};<4.5mm,5mm>*{}**@{-},
 <0.6mm,0.44mm>*{};<8mm,5mm>*{}**@{-},
   <0mm,0mm>*{};<-10.5mm,5.9mm>*{^{\sigma(1)}}**@{},
   <0mm,0mm>*{};<-4mm,5.9mm>*{^{\sigma(2)}}**@{},
   <0mm,0mm>*{};<10.0mm,5.9mm>*{^{\sigma(m)}}**@{},
 <-0.6mm,-0.44mm>*{};<-8mm,-5mm>*{}**@{-},
 <-0.4mm,-0.7mm>*{};<-4.5mm,-5mm>*{}**@{-},
 <0mm,0mm>*{};<1mm,-5mm>*{\ldots}**@{},
 <0.4mm,-0.7mm>*{};<4.5mm,-5mm>*{}**@{-},
 <0.6mm,-0.44mm>*{};<8mm,-5mm>*{}**@{-},
   <0mm,0mm>*{};<-10.5mm,-6.9mm>*{^{\tau(1)}}**@{},
   <0mm,0mm>*{};<-4mm,-6.9mm>*{^{\tau(2)}}**@{},
   <0mm,0mm>*{};<10.0mm,-6.9mm>*{^{\tau(n)}}**@{},
 \end{xy}}\Ea
=(-1)^{c|\sigma|+d|\tau|}
\Ba{c}\resizebox{14mm}{!}{\begin{xy}
 <0mm,0mm>*{\circ};<0mm,0mm>*{}**@{},
 <-0.6mm,0.44mm>*{};<-8mm,5mm>*{}**@{-},
 <-0.4mm,0.7mm>*{};<-4.5mm,5mm>*{}**@{-},
 <0mm,0mm>*{};<-1mm,5mm>*{\ldots}**@{},
 <0.4mm,0.7mm>*{};<4.5mm,5mm>*{}**@{-},
 <0.6mm,0.44mm>*{};<8mm,5mm>*{}**@{-},
   <0mm,0mm>*{};<-8.5mm,5.5mm>*{^1}**@{},
   <0mm,0mm>*{};<-5mm,5.5mm>*{^2}**@{},
   <0mm,0mm>*{};<4.5mm,5.5mm>*{^{m\hspace{-0.5mm}-\hspace{-0.5mm}1}}**@{},
   <0mm,0mm>*{};<9.0mm,5.5mm>*{^m}**@{},
 <-0.6mm,-0.44mm>*{};<-8mm,-5mm>*{}**@{-},
 <-0.4mm,-0.7mm>*{};<-4.5mm,-5mm>*{}**@{-},
 <0mm,0mm>*{};<-1mm,-5mm>*{\ldots}**@{},
 <0.4mm,-0.7mm>*{};<4.5mm,-5mm>*{}**@{-},
 <0.6mm,-0.44mm>*{};<8mm,-5mm>*{}**@{-},
   <0mm,0mm>*{};<-8.5mm,-6.9mm>*{^1}**@{},
   <0mm,0mm>*{};<-5mm,-6.9mm>*{^2}**@{},
   <0mm,0mm>*{};<4.5mm,-6.9mm>*{^{n\hspace{-0.5mm}-\hspace{-0.5mm}1}}**@{},
   <0mm,0mm>*{};<9.0mm,-6.9mm>*{^n}**@{},
 \end{xy}}\Ea \ \ \forall \sigma\in \bS_m, \forall\tau\in \bS_n,\ \ m,n\geq 1, m+n\geq 3,
\Eeq
with $m+n\geq 3$,$m,n\geq 1$,
and equipped with the differential
given on the generators by
\Beq\label{differential in HoLBcd}
\delta
\Ba{c}\resizebox{14mm}{!}{\begin{xy}
 <0mm,0mm>*{\circ};<0mm,0mm>*{}**@{},
 <-0.6mm,0.44mm>*{};<-8mm,5mm>*{}**@{-},
 <-0.4mm,0.7mm>*{};<-4.5mm,5mm>*{}**@{-},
 <0mm,0mm>*{};<-1mm,5mm>*{\ldots}**@{},
 <0.4mm,0.7mm>*{};<4.5mm,5mm>*{}**@{-},
 <0.6mm,0.44mm>*{};<8mm,5mm>*{}**@{-},
   <0mm,0mm>*{};<-8.5mm,5.5mm>*{^1}**@{},
   <0mm,0mm>*{};<-5mm,5.5mm>*{^2}**@{},
   <0mm,0mm>*{};<4.5mm,5.5mm>*{^{m\hspace{-0.5mm}-\hspace{-0.5mm}1}}**@{},
   <0mm,0mm>*{};<9.0mm,5.5mm>*{^m}**@{},
 <-0.6mm,-0.44
 mm>*{};<-8mm,-5mm>*{}**@{-},
 <-0.4mm,-0.7mm>*{};<-4.5mm,-5mm>*{}**@{-},
 <0mm,0mm>*{};<-1mm,-5mm>*{\ldots}**@{},
 <0.4mm,-0.7mm>*{};<4.5mm,-5mm>*{}**@{-},
 <0.6mm,-0.44mm>*{};<8mm,-5mm>*{}**@{-},
   <0mm,0mm>*{};<-8.5mm,-6.9mm>*{^1}**@{},
   <0mm,0mm>*{};<-5mm,-6.9mm>*{^2}**@{},
   <0mm,0mm>*{};<4.5mm,-6.9mm>*{^{n\hspace{-0.5mm}-\hspace{-0.5mm}1}}**@{},
   <0mm,0mm>*{};<9.0mm,-6.9mm>*{^n}**@{},
 \end{xy}}\Ea
\ \ = \ \
 \sum_{[1,\ldots,m]=I_1\sqcup I_2\atop
 {|I_1|\geq 0, |I_2|\geq 1}}
 \sum_{[1,\ldots,n]=J_1\sqcup J_2\atop
 {|J_1|\geq 1, |J_2|\geq 1}
}\hspace{0mm}
(-1)^{d(\# J_1 + \# I_1\# J_2 + \sigma(I_1,I_2) + \sigma(J_1,J_2))}
\Ba{c}\resizebox{22mm}{!}{ \begin{xy}
 <0mm,0mm>*{\circ};<0mm,0mm>*{}**@{},
 <-0.6mm,0.44mm>*{};<-8mm,5mm>*{}**@{-},
 <-0.4mm,0.7mm>*{};<-4.5mm,5mm>*{}**@{-},
 <0mm,0mm>*{};<0mm,5mm>*{\ldots}**@{},
 <0.4mm,0.7mm>*{};<4.5mm,5mm>*{}**@{-},
 <0.6mm,0.44mm>*{};<12.4mm,4.8mm>*{}**@{-},
     <0mm,0mm>*{};<-2mm,7mm>*{\overbrace{\ \ \ \ \ \ \ \ \ \ \ \ }}**@{},
     <0mm,0mm>*{};<-2mm,9mm>*{^{I_1}}**@{},
 <-0.6mm,-0.44mm>*{};<-8mm,-5mm>*{}**@{-},
 <-0.4mm,-0.7mm>*{};<-4.5mm,-5mm>*{}**@{-},
 <0mm,0mm>*{};<-1mm,-5mm>*{\ldots}**@{},
 <0.4mm,-0.7mm>*{};<4.5mm,-5mm>*{}**@{-},
 <0.6mm,-0.44mm>*{};<8mm,-5mm>*{}**@{-},
      <0mm,0mm>*{};<0mm,-7mm>*{\underbrace{\ \ \ \ \ \ \ \ \ \ \ \ \ \ \
      }}**@{},
      <0mm,0mm>*{};<0mm,-10.6mm>*{_{J_1}}**@{},
 <13mm,5mm>*{};<13mm,5mm>*{\circ}**@{},
 <12.6mm,5.44mm>*{};<5mm,10mm>*{}**@{-},
 <12.6mm,5.7mm>*{};<8.5mm,10mm>*{}**@{-},
 <13mm,5mm>*{};<13mm,10mm>*{\ldots}**@{},
 <13.4mm,5.7mm>*{};<16.5mm,10mm>*{}**@{-},
 <13.6mm,5.44mm>*{};<20mm,10mm>*{}**@{-},
      <13mm,5mm>*{};<13mm,12mm>*{\overbrace{\ \ \ \ \ \ \ \ \ \ \ \ \ \ }}**@{},
      <13mm,5mm>*{};<13mm,14mm>*{^{I_2}}**@{},
 <12.4mm,4.3mm>*{};<8mm,0mm>*{}**@{-},
 <12.6mm,4.3mm>*{};<12mm,0mm>*{\ldots}**@{},
 <13.4mm,4.5mm>*{};<16.5mm,0mm>*{}**@{-},
 <13.6mm,4.8mm>*{};<20mm,0mm>*{}**@{-},
     <13mm,5mm>*{};<14.3mm,-2mm>*{\underbrace{\ \ \ \ \ \ \ \ \ \ \ }}**@{},
     <13mm,5mm>*{};<14.3mm,-4.5mm>*{_{J_2}}**@{},
 \end{xy}}\Ea
\Eeq
where the vertices on the r.h.s.\ are ordered in such a way that the lowest one comes first.
A representation $\rho$ of $\HoLBcd$ in a dg vector space $V$ can be identified with
a degree $c+d+1$ Maurer-Cartan element $\pi$, $\{\pi,\pi\}=0$,
 in the completed graded commutative algebra
$$
\pi=\sum_{m,n\geq 1} \pi_{n}^{m} \in \prod_{m,n\geq 1,m+n\geq 3}
\Hom(\odot^n (V[d]), \odot^m (V[-c])\subset \prod_{k\geq 0}\odot^k (V^*[-d]\oplus V[-c])
$$
equipped with the Poisson type Lie bracket $\{\ ,\ \}$, of degree $-c-d$ induced by the natural paring
$V\ot V^*\rar \K$. Here, by definition,
 $$
 \pi_{n}^{m}:=\rho\left(\Ba{c}\resizebox{14mm}{!}{\begin{xy}
 <0mm,0mm>*{\circ};<0mm,0mm>*{}**@{},
 <-0.6mm,0.44mm>*{};<-8mm,5mm>*{}**@{-},
 <-0.4mm,0.7mm>*{};<-4.5mm,5mm>*{}**@{-},
 <0mm,0mm>*{};<-1mm,5mm>*{\ldots}**@{},
 <0.4mm,0.7mm>*{};<4.5mm,5mm>*{}**@{-},
 <0.6mm,0.44mm>*{};<8mm,5mm>*{}**@{-},
   <0mm,0mm>*{};<-8.5mm,5.5mm>*{^1}**@{},
   <0mm,0mm>*{};<-5mm,5.5mm>*{^2}**@{},
   <0mm,0mm>*{};<4.5mm,5.5mm>*{^{m\hspace{-0.5mm}-\hspace{-0.5mm}1}}**@{},
   <0mm,0mm>*{};<9.0mm,5.5mm>*{^m}**@{},
 <-0.6mm,-0.44
 mm>*{};<-8mm,-5mm>*{}**@{-},
 <-0.4mm,-0.7mm>*{};<-4.5mm,-5mm>*{}**@{-},
 <0mm,0mm>*{};<-1mm,-5mm>*{\ldots}**@{},
 <0.4mm,-0.7mm>*{};<4.5mm,-5mm>*{}**@{-},
 <0.6mm,-0.44mm>*{};<8mm,-5mm>*{}**@{-},
   <0mm,0mm>*{};<-8.5mm,-6.9mm>*{^1}**@{},
   <0mm,0mm>*{};<-5mm,-6.9mm>*{^2}**@{},
   <0mm,0mm>*{};<4.5mm,-6.9mm>*{^{n\hspace{-0.5mm}-\hspace{-0.5mm}1}}**@{},
   <0mm,0mm>*{};<9.0mm,-6.9mm>*{^n}**@{},
 \end{xy}}\Ea  \right)\ \ \ \text{for}\ m+n\geq 3,
 $$
 and $\pi_1^1$ is the given differential in $V$.
 We are mostly interested in this paper in the case $c=1$, $d=0$, $V=\R^N$ when $\pi$ becomes a formal Poisson structure on $\R^N$ which vanishes at $0\in \R^N$,
\Beq\label{2: pi Poisson}
\pi=\sum_{a,b=1}^N\pi^{a,b} (x)\p_a\wedge \p_b, \ \  \text{with}\ \ \pi^{ab} (x) :=\sum_{n\geq 1}\frac{1}{n!}\pi^{ab}_{c_1\ldots c_n} x^{c_1}\cdots x^{c_n}
\Eeq
for some constants $\pi^{ab}_{c_1\ldots c_n}=-\pi^{ba}_{c_1\ldots c_n}\in \R$.

\sip

Consider a differential ideal $I$  in $\HoLB_{c,d}$ generated by graphs which
 contain at least one $(m,n)$-corolla with $m\neq n$, and let
$$
\HoQP_{c,d}:=\HoLBcd/I
$$
be the quotient dg properad. It is a free properad generated by $(m,m)$-corollas with $m\geq 2$ from the family (\ref{2: symmetries of HoLiebcd corollas}),
$$
\Ba{c}\resizebox{17mm}{!}{\begin{xy}
 <0mm,0mm>*{\circ};<0mm,0mm>*{}**@{},
 <-0.6mm,0.44mm>*{};<-8mm,5mm>*{}**@{-},
 <-0.4mm,0.7mm>*{};<-4.5mm,5mm>*{}**@{-},
 <0mm,0mm>*{};<-1mm,5mm>*{\ldots}**@{},
 <0.4mm,0.7mm>*{};<4.5mm,5mm>*{}**@{-},
 <0.6mm,0.44mm>*{};<8mm,5mm>*{}**@{-},
   <0mm,0mm>*{};<-8.5mm,5.5mm>*{^1}**@{},
   <0mm,0mm>*{};<-5mm,5.5mm>*{^2}**@{},
   <0mm,0mm>*{};<4.5mm,5.5mm>*{^{m\hspace{-0.5mm}-\hspace{-0.5mm}1}}**@{},
   <0mm,0mm>*{};<9.0mm,5.5mm>*{^m}**@{},
 <-0.6mm,-0.44
 mm>*{};<-8mm,-5mm>*{}**@{-},
 <-0.4mm,-0.7mm>*{};<-4.5mm,-5mm>*{}**@{-},
 <0mm,0mm>*{};<-1mm,-5mm>*{\ldots}**@{},
 <0.4mm,-0.7mm>*{};<4.5mm,-5mm>*{}**@{-},
 <0.6mm,-0.44mm>*{};<8mm,-5mm>*{}**@{-},
   <0mm,0mm>*{};<-8.5mm,-6.9mm>*{^1}**@{},
   <0mm,0mm>*{};<-5mm,-6.9mm>*{^2}**@{},
   <0mm,0mm>*{};<4.5mm,-6.9mm>*{^{m\hspace{-0.5mm}-\hspace{-0.5mm}1}}**@{},
   <0mm,0mm>*{};<9.0mm,-6.9mm>*{^m}**@{},
 \end{xy}}\Ea, \ \ \ \ m\geq 2
$$
and equipped with the  differential $\delta$ induced from (\ref{differential in HoLBcd}); it is easy to see that $\delta$ acts trivially on the $(2,2)$-corolla. Let $J$ be the differential
closure of the ideal in $\HoQP_{c,d}$ generated by $(m,m)$-corollas with $m\geq 3$ and set
$$
\QP_{c,d}:=\HoQP_{c,d}/J
$$
The latter properad is a quadratic properad generated by the $(2,2)$-corolla of degree $1-c-d$,
$$
\Ba{c}\begin{xy}
 <0.5mm,0.5mm>*{};<2.2mm,2.2mm>*{}**@{-},
 <-0.48mm,0.48mm>*{};<-2.2mm,2.2mm>*{}**@{-},
 <0.5mm,-0.5mm>*{};<2.2mm,-2.2mm>*{}**@{-},
 <-0.48mm,-0.48mm>*{};<-2.2mm,-2.2mm>*{}**@{-},
 <0mm,0mm>*{\circ};
 <2.7mm,2.8mm>*{^{_2}}**@{},
 <-2.7mm,2.8mm>*{^{_1}}**@{},
 <2.8mm,-3.3mm>*{_{_2}}**@{},
 <-2.7mm,-3.3mm>*{_{_1}}**@{},
 \end{xy}\Ea
 =(-1)^{c}
 \Ba{c}\begin{xy}
 <0.5mm,0.5mm>*{};<2.2mm,2.2mm>*{}**@{-},
 <-0.48mm,0.48mm>*{};<-2.2mm,2.2mm>*{}**@{-},
 <0.5mm,-0.5mm>*{};<2.2mm,-2.2mm>*{}**@{-},
 <-0.48mm,-0.48mm>*{};<-2.2mm,-2.2mm>*{}**@{-},
 <0mm,0mm>*{\circ};
 <2.7mm,2.8mm>*{^{_1}}**@{},
 <-2.7mm,2.8mm>*{^{_2}}**@{},
 <2.8mm,-3.3mm>*{_{_2}}**@{},
 <-2.7mm,-3.3mm>*{_{_1}}**@{},
 \end{xy}\Ea
 =(-1)^d
 \Ba{c}\begin{xy}
 <0.5mm,0.5mm>*{};<2.2mm,2.2mm>*{}**@{-},
 <-0.48mm,0.48mm>*{};<-2.2mm,2.2mm>*{}**@{-},
 <0.5mm,-0.5mm>*{};<2.2mm,-2.2mm>*{}**@{-},
 <-0.48mm,-0.48mm>*{};<-2.2mm,-2.2mm>*{}**@{-},
 <0mm,0mm>*{\circ};
 <2.7mm,2.8mm>*{^{_2}}**@{},
 <-2.7mm,2.8mm>*{^{_1}}**@{},
 <2.8mm,-3.3mm>*{_{_1}}**@{},
 <-2.7mm,-3.3mm>*{_{_2}}**@{},
 \end{xy}\Ea
$$
modulo the relations
$$
\sum_{\tau,\sigma\in \Z_3}
\Ba{c}
{ \begin{xy}
<0.5mm,0.5mm>*{};<2.4mm,2.4mm>*{}**@{-},
 <-0.48mm,0.48mm>*{};<-2.4mm,2.4mm>*{}**@{-},
 <0mm,0mm>*{\circ};
 <0.39mm,-0.39mm>*{};<2.4mm,-2.4mm>*{}**@{-},
 <-0.35mm,-0.35mm>*{};<-1.9mm,-1.9mm>*{}**@{-},
 <5.0mm,3.5mm>*{^{_{\tau(3)}}}**@{},
<-2.5mm,3.5mm>*{^{_{\tau(2)}}}**@{},
<-5.6mm,1.0mm>*{^{_{\tau(1)}}}**@{},
 <-2.4mm,-2.4mm>*{\circ};
 <-2.8mm,-2.0mm>*{};<-4.8mm,0.0mm>*{}**@{-},
 <-2.0mm,-2.8mm>*{};<0mm,-4.9mm>*{}**@{-},
 <-2.8mm,-2.9mm>*{};<-4.7mm,-4.9mm>*{}**@{-},
<5.0mm,-4.0mm>*{^{_{\sigma(3)}}}**@{},
<1.5mm,-6.7mm>*{^{_{\sigma(2)}}}**@{},
<-5.2mm,-6.7mm>*{^{_{\sigma(1)}}}**@{},
 \end{xy}}\Ea=0
$$
where $\Z_3$ is the subgroup of $\bS_3$ generated by the permutation $(123)$. Thus
 representations of the properad $\QP_{1,0}$ in an arbitrary (not necessarily finite-dimensional) graded vector space $V$ can be identified with quadratic Poisson structures on the vector space $V^*$; in particular, in the case $V=\R^N$ they are given by (\ref{2: pi Poisson}) with only $\pi_{c_1c_2}^{ab}\neq 0$. It is shown in \cite{KM} that the canonical epimorphism
$$
p: \HoQP_{c,d} \lon \QP_{c,d}
$$
is a quasi-isomorphism implying that the properad $\QP_{c,d}$ is Koszul (we do not use this reult in this paper).

\sip

\sip

The category of wheeled props has been introduced  and studied in \cite{Me1, MMS}. The properads $\HoQP_{c,d}$ and $\QP_{c,d}$ have very simple wheeled closures denoted by $\wHoQP_{c,d}$ and $\QP_{c,d}^\circlearrowright$ respectively. They have the same sets of generators, but one is a allowed to build from them (via gluing outgoing legs of one generator to ingoing legs of another or the same generator) graphs which have closed paths of directed edges. For, example the graph
$$
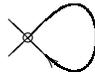
\begin{xy}
 <0mm,0mm>*{\circ};<0mm,0mm>*{}**@{},
<0mm,0mm>*{};<2.5mm,-2.5mm>*{}**@{-},
<0mm,0mm>*{};<-2.5mm,-2.5mm>*{}**@{-},
<0mm,0mm>*{};<2.5mm,2.5mm>*{}**@{-},
<0mm,0mm>*{};<-2.5mm,2.5mm>*{}**@{-},
<0mm,0mm>*{};<2.5mm,2.5mm>*{}**@{-},
(2.4,2.3)*{}
   \ar@{->}@(ur,dr) (2.4,-2.3)*{}
\end{xy}
$$
is allowed in both properads $\wHoQP_{c,d}$ and $\QP_{c,d}^\circlearrowright$, and, in the notation used in (\ref{2: pi Poisson}), gets represented by the linear vector field
$$
\sum_{a,b,c} \frac{1}{2} \pi_{bc}^{ac}x^b \p_a.
$$
 This example shows a fundamental difference between ordinary props and their wheeled closures: if, for example, the properad $\QP_{0,1}$ controls quadratic Poisson structures in both finite- and infinite-dimensional spaces, its wheeled closure $\QP_{c,d}^\circlearrowright$ admits representations only in {\em finite}-dimensional vector spaces as it involves the trace operation $\Hom(V,V)\rar \K$! Hence in the context of deformation quantization we should be interested in properads $\wHoQP_{0,1}$ and $\QP_{0,1}^\circlearrowright$ as universal quantization formulae {\em must}\, use graphs with wheels (see \S 5 in \cite{Me3}).

\sip

There is an exact functor from the category of properads to the category of props \cite{V}.  Properads $\cP$ studied in this paper are spanned by {\em connected}\, graphs while their prop closures use disjoint unions of connected graphs, i.e.\ they have the same generators and relations but the graphs built from generators are not necessarily connected.
The prop enveloping of a properad $\cP$ is denoted by the same symbol $\cP$; often it plays no role whether we work with a properad $\cP$ or its prop enveloping $\cP$ but when it is important we say this explicitly.

\sip

Another important technical point is that deformation quantization formulae \cite{Ko2} of a generic $\Z$-graded (quadratic) Poisson structure involve formal power series in an auxiliary formal parameter $\hbar$ (the ``Planck constant"). In terms of props this fact leads us to consider {\em genus completed}\, version, $\wiweHoQP_{0,1}$ of the prop $\weHoQP_{0,1}$. This is one of the main reasons why we develop below deformation theory of completed  wheeled props  $\wiweHoQP_{c,d}$  rather than of the ordinary ones.

\subsection{Complexes of derivations of properads}  There is a useful endofunctor,
$\cP \rar \cP^+$, in the category of (wheeled) dg props introduced in \cite{Me5}. Given any dg (wheeled) prop $\cP$, the associated extended dg prop
$\cP^+$ is obtained from $\cP$ by adjoining an extra generator $t\in \cP^{+}(1,1)$
 and is uniquely characterized by the property: there is a 1-1 correspondence between representations
 $$
 \rho: \cP^+ \lon \cE nd_V
 $$
of $\cP^+$ in a dg vector space $(V,d)$, and representations of $\cP$ in the same space $V$
but equipped with a deformed differential $d+d'$, where $d':=\rho(t)$.
The complex of derivations of an arbitrary prop $\cP$ has been defined in \cite{MW1}  as the space of derivations,
$$
\Der(\cP):= \Der(\cP^+ \rar \cP)
$$
of the plus-extended prop with values in $\cP$; the plus extension is used in order not to loose an important information about the homotopy theory of the (wheeled) prop $\cP$, cf.\ \cite{MW1,MW2,GY,MW3}. The dg props $\widehat{\HoQP}_{c,d}^{+}$ an ${\HoQP}_{c,d}^{+\circlearrowright}$ is easy to describe explicitly --- both are built from generators (\ref{2: symmetries of HoLiebcd corollas}) and one extra
$(1,1)$-generator
$
\begin{xy}
 <0mm,-0.55mm>*{};<0mm,-3mm>*{}**@{-},
 <0mm,0.5mm>*{};<0mm,3mm>*{}**@{-},
 <0mm,0mm>*{\bu};<0mm,0mm>*{}**@{},
 \end{xy}
 $ which is assigned degree $+1$.
The differential $\delta^+$ is given on the new generator by the formula
 $$
 \delta^+
 \begin{xy}
 <0mm,-0.55mm>*{};<0mm,-3mm>*{}**@{-},
 <0mm,0.5mm>*{};<0mm,3mm>*{}**@{-},
 <0mm,0mm>*{\bu};<0mm,0mm>*{}**@{},
 \end{xy} :=
 \Ba{c}
 \begin{xy}
 <0mm,-0.5mm>*{};<0mm,-3.5mm>*{}**@{-},
 <0mm,0.5mm>*{};<0mm,3.5mm>*{}**@{-},
 <0mm,4.5mm>*{};<0mm,7mm>*{}**@{-},
 <0mm,0mm>*{\bu};
 <0mm,4mm>*{\bu};
 \end{xy}
 \Ea
 $$
 while its value on ``old" generators is given by
 $$
 \delta^+
 \Ba{c}\resizebox{18mm}{!}{\begin{xy}
 <0mm,0mm>*{\circ};<0mm,0mm>*{}**@{},
 <-0.6mm,0.44mm>*{};<-8mm,5mm>*{}**@{-},
 <-0.4mm,0.7mm>*{};<-4.5mm,5mm>*{}**@{-},
 <0mm,0mm>*{};<-1mm,5mm>*{\ldots}**@{},
 <0.4mm,0.7mm>*{};<4.5mm,5mm>*{}**@{-},
 <0.6mm,0.44mm>*{};<8mm,5mm>*{}**@{-},
   <0mm,0mm>*{};<-8.5mm,5.5mm>*{^1}**@{},
   <0mm,0mm>*{};<-5mm,5.5mm>*{^2}**@{},
   <0mm,0mm>*{};<4.5mm,5.5mm>*{^{m\hspace{-0.5mm}-\hspace{-0.5mm}1}}**@{},
   <0mm,0mm>*{};<9.0mm,5.5mm>*{^m}**@{},
 <-0.6mm,-0.44
 mm>*{};<-8mm,-5mm>*{}**@{-},
 <-0.4mm,-0.7mm>*{};<-4.5mm,-5mm>*{}**@{-},
 <0mm,0mm>*{};<-1mm,-5mm>*{\ldots}**@{},
 <0.4mm,-0.7mm>*{};<4.5mm,-5mm>*{}**@{-},
 <0.6mm,-0.44mm>*{};<8mm,-5mm>*{}**@{-},
   <0mm,0mm>*{};<-8.5mm,-6.9mm>*{^1}**@{},
   <0mm,0mm>*{};<-5mm,-6.9mm>*{^2}**@{},
   <0mm,0mm>*{};<4.5mm,-6.9mm>*{^{m\hspace{-0.5mm}-\hspace{-0.5mm}1}}**@{},
   <0mm,0mm>*{};<9.0mm,-6.9mm>*{^m}**@{},
 \end{xy}}\Ea
:= \delta
 \Ba{c}\resizebox{18mm}{!}{\begin{xy}
 <0mm,0mm>*{\circ};<0mm,0mm>*{}**@{},
 <-0.6mm,0.44mm>*{};<-8mm,5mm>*{}**@{-},
 <-0.4mm,0.7mm>*{};<-4.5mm,5mm>*{}**@{-},
 <0mm,0mm>*{};<-1mm,5mm>*{\ldots}**@{},
 <0.4mm,0.7mm>*{};<4.5mm,5mm>*{}**@{-},
 <0.6mm,0.44mm>*{};<8mm,5mm>*{}**@{-},
   <0mm,0mm>*{};<-8.5mm,5.5mm>*{^1}**@{},
   <0mm,0mm>*{};<-5mm,5.5mm>*{^2}**@{},
   <0mm,0mm>*{};<4.5mm,5.5mm>*{^{m\hspace{-0.5mm}-\hspace{-0.5mm}1}}**@{},
   <0mm,0mm>*{};<9.0mm,5.5mm>*{^m}**@{},
 <-0.6mm,-0.44
 mm>*{};<-8mm,-5mm>*{}**@{-},
 <-0.4mm,-0.7mm>*{};<-4.5mm,-5mm>*{}**@{-},
 <0mm,0mm>*{};<-1mm,-5mm>*{\ldots}**@{},
 <0.4mm,-0.7mm>*{};<4.5mm,-5mm>*{}**@{-},
 <0.6mm,-0.44mm>*{};<8mm,-5mm>*{}**@{-},
   <0mm,0mm>*{};<-8.5mm,-6.9mm>*{^1}**@{},
   <0mm,0mm>*{};<-5mm,-6.9mm>*{^2}**@{},
   <0mm,0mm>*{};<4.5mm,-6.9mm>*{^{m\hspace{-0.5mm}-\hspace{-0.5mm}1}}**@{},
   <0mm,0mm>*{};<9.0mm,-6.9mm>*{^m}**@{},
 \end{xy}}\Ea
+
\overset{m-1}{\underset{i=0}{\sum}} \ \
 \resizebox{18mm}{!}{\begin{xy}
 <0mm,0mm>*{\circ};<0mm,0mm>*{}**@{},
 <-0.5mm,0.5mm>*{};<-8mm,5mm>*{}**@{-},
 <-0.5mm,0.5mm>*{};<-3.5mm,5mm>*{}**@{-},
 <0.5mm,0.4mm>*{};<-6mm,5mm>*{..}**@{},
 <0mm,0.55mm>*{};<0mm,9mm>*{}**@{-},
  <0mm,5mm>*{\bullet};
  <0mm,5mm>*{};<0mm,10mm>*{^{i\hspace{-0.2mm}+\hspace{-0.5mm}1}}**@{},
<0mm,0.5mm>*{};<8mm,5mm>*{}**@{-},
<0.4mm,0.4mm>*{};<3.5mm,5mm>*{}**@{-},
 <0mm,0mm>*{};<6mm,5mm>*{..}**@{},
   <0mm,0mm>*{};<-8.5mm,5.5mm>*{^1}**@{},
   <0mm,0mm>*{};<-4mm,5.5mm>*{^i}**@{},
   <0mm,0mm>*{};<9.0mm,5.5mm>*{^m}**@{},
 <-0.5mm,-0.5mm>*{};<-8mm,-5mm>*{}**@{-},
 <-0.5mm,-0.5mm>*{};<-4.5mm,-5mm>*{}**@{-},
 <0mm,0mm>*{};<-1mm,-5mm>*{\ldots}**@{},
 <0.5mm,-0.5mm>*{};<4.5mm,-5mm>*{}**@{-},
 <0.5mm,-0.5mm>*{};<8mm,-5mm>*{}**@{-},
   <0mm,0mm>*{};<-8.5mm,-6.9mm>*{^1}**@{},
   <0mm,0mm>*{};<-5mm,-6.9mm>*{^2}**@{},
   <0mm,0mm>*{};<4.5mm,-6.9mm>*{^{m\hspace{-0.5mm}-\hspace{-0.5mm}1}}**@{},
   <0mm,0mm>*{};<9.0mm,-6.9mm>*{^m}**@{},
 \end{xy}
 }
 \ \
 \pm
 \ \
\overset{m-1}{\underset{i=0}{\sum}}\ \
\resizebox{18mm}{!}{\begin{xy}
 <0mm,0mm>*{\circ};<0mm,0mm>*{}**@{},
 <-0.5mm,-0.5mm>*{};<-8mm,-5mm>*{}**@{-},
 <-0.5mm,-0.5mm>*{};<-3.5mm,-5mm>*{}**@{-},
 <0.5mm,-0.4mm>*{};<-6mm,-5mm>*{..}**@{},
 <0mm,-0.55mm>*{};<0mm,-9mm>*{}**@{-},
  <0mm,-5mm>*{\bullet};
  <0mm,-5mm>*{};<0mm,-10.9mm>*{^{i\hspace{-0.2mm}+\hspace{-0.5mm}1}}**@{},
<0mm,-0.5mm>*{};<8mm,-5mm>*{}**@{-},
<0.4mm,-0.4mm>*{};<3.5mm,-5mm>*{}**@{-},
 <0mm,0mm>*{};<6mm,-5mm>*{..}**@{},
   <0mm,0mm>*{};<-8.5mm,-5.9mm>*{_1}**@{},
   <0mm,0mm>*{};<-4mm,-5.9mm>*{_i}**@{},
   <0mm,0mm>*{};<9.0mm,-5.9mm>*{_m}**@{},
 <-0.5mm,0.5mm>*{};<-8mm,5mm>*{}**@{-},
 <-0.5mm,0.5mm>*{};<-4.5mm,5mm>*{}**@{-},
 <0mm,0mm>*{};<-1mm,5mm>*{\ldots}**@{},
 <0.5mm,0.5mm>*{};<4.5mm,5mm>*{}**@{-},
 <0.5mm,0.5mm>*{};<8mm,5mm>*{}**@{-},
   <0mm,0mm>*{};<-8.5mm,6.9mm>*{_1}**@{},
   <0mm,0mm>*{};<-5mm,6.9mm>*{_2}**@{},
   <0mm,0mm>*{};<4.5mm,6.9mm>*{_{m\hspace{-0.5mm}-\hspace{-0.5mm}1}}**@{},
   <0mm,0mm>*{};<9.0mm,6.9mm>*{_m}**@{},
 \end{xy}
 }
 $$

It is not hard to show that the complexes $\Der(\cP^+ \rar \cP)$ and
$\Der(\cP^+ \rar \cP^+)$ are quasi-isomorphic so that there is no loss of generality to define a derivation complex of a dg prop $\cP$ as
$\Der(\cP^+ \rar \cP^+)$
making the Lie algebra structure on it more transparent. Hence it is a matter of taste
which definition to use. We prefer the original one from (\cite{MW1}) and set
$$
\Der(\HoQP_{c,d}) := \Der(\widehat{\HoQP}_{c,d}^{\atop+} \rar \widehat{\HoQP}_{c,d}^{\atop}), \ \ \
\Der(\weHoQPcd) := \Der(\widehat{\HoQP}_{c,d}^{\atop +\circlearrowright}
\rar \widehat{\HoQP}_{c,d}^{\atop \circlearrowright}).
$$
As graded vector spaces, they can be identified with spaces generated by certain graphs. For example
\Beq\label{2: Der(weHoQP) as vector space}
\Der(\weHoQPcd)=\prod_{m\geq 1} \left(\widehat{\HoQP}_{c,d}^{\atop \circlearrowright}(m,m) \otimes \sgn_m^{\ot |c|}\otimes \sgn_m^{\ot |d|}\right)^{\bS_m\times \bS_m}[1+(c+d)(1-m)],
\Eeq
Its elements are directed not necessary connected graphs (possibly, with wheels) which might  have incoming and outgoing legs and vertices have the same numbers of incoming and outgoing edges, for example
\Beq\label{2: pic for graphs in Der Hoqios}
\Ba{c}
\resizebox{15mm}{!}{ \xy
(0,0)*{\bu}="d1",
(10,0)*{\bu}="d2",
(15,3)*{}="d2c",
(-5,-5)*{}="d1l",
(-5,-5)*{}="dl",
(5,-5)*{}="dc",
(15,-5)*{}="dr",
(0,10)*{\bu}="u1",
(10,10)*{\bu}="u2",
(15,6)*{}="u2r",
(5,15)*{}="uc",
(15,15)*{}="ur",
(0,15)*{}="ul",
\ar @{<-} "d1";"d2" <0pt>
\ar @{<-} "d1";"d1l" <0pt>
\ar @{<-} "d2";"dc" <0pt>
\ar @{->} "d2";"d2c" <0pt>
\ar @{<-} "d2";"dr" <0pt>
\ar @{<-} "u1";"d1" <0pt>
\ar @{->} "u1";"u2" <0pt>
\ar @{<-} "u1";"d2" <0pt>
\ar @{->} "u2";"d2" <0pt>
\ar @{<-} "u2";"d1" <0pt>
\ar @{<-} "u2";"u2r" <0pt>
\ar @{<-} "uc";"u2" <0pt>
\ar @{<-} "ur";"u2" <0pt>
\ar @{<-} "ul";"u1" <0pt>
\endxy}
\Ea
 \in \Der(\weHoQPcd)
\Eeq
Note that
contrary to the prop $\HoLBcd^\circlearrowright$, the out- or ingoing legs (if any) of a graph from $\Der(\weHoQPcd)$ have their numerical labels are (skew)symmetrized in accordance with the parity of the integer parameters $c$ and $d$; one can think that they are not labelled at all --- just for odd $c$ or $d$  some ordering is chosen (up to an even permutation). The subspace $\Der(\HoQP_{c,d})\subset \Der(\weHoQPcd)$ is generated by similar graphs but with no wheels (in this case the summand above corresponding to $m=0$ vanishes).

\sip

 The value of the differential $d$ on an element $\Ga\in \Der(\weHoQPcd)$ consists of three terms
\Beq\label{d in Der(Holieb)}
d \Gamma =
 \delta\Gamma
  \pm
  \sum_{m\geq 2}
  \Ba{c}\resizebox{14mm}{!}{\begin{xy}
 <0mm,0mm>*{\bu};<0mm,0mm>*{}**@{},
 <-0.6mm,0.44mm>*{};<-8mm,5mm>*{}**@{-},
 <-0.4mm,0.7mm>*{};<-4.5mm,5mm>*{}**@{-},
 <0mm,0mm>*{};<-1mm,5mm>*{\ldots}**@{},
 <0.4mm,0.7mm>*{};<4.5mm,5mm>*{}**@{-},
 <0.6mm,0.44mm>*{};<10mm,6mm>*{}**@{-},
   <0mm,0mm>*{};<12.0mm,7.5mm>*{\Ga}**@{},
 <-0.6mm,-0.44mm>*{};<-8mm,-5mm>*{}**@{-},
 <-0.4mm,-0.7mm>*{};<-4.5mm,-5mm>*{}**@{-},
 <0mm,0mm>*{};<-1mm,-5mm>*{\ldots}**@{},
 <0.4mm,-0.7mm>*{};<4.5mm,-5mm>*{}**@{-},
 <0.6mm,-0.44mm>*{};<8mm,-5mm>*{}**@{-},
 \end{xy}}\Ea
   \pm
  \sum
  \Ba{c}\resizebox{14mm}{!}{\begin{xy}
 <0mm,0mm>*{\bu};<0mm,0mm>*{}**@{},
 <-0.6mm,0.44mm>*{};<-8mm,5mm>*{}**@{-},
 <-0.4mm,0.7mm>*{};<-4.5mm,5mm>*{}**@{-},
 <0mm,0mm>*{};<-1mm,5mm>*{\ldots}**@{},
 <0.4mm,0.7mm>*{};<4.5mm,5mm>*{}**@{-},
 <0.6mm,0.44mm>*{};<-10mm,-6mm>*{}**@{-},
   <0mm,0mm>*{};<-12.0mm,-7.5mm>*{\Ga}**@{},
 <-0.6mm,-0.44mm>*{};<8mm,5mm>*{}**@{-},
 <-0.4mm,-0.7mm>*{};<-4.5mm,-5mm>*{}**@{-},
 <0mm,0mm>*{};<-1mm,-5mm>*{\ldots}**@{},
 <0.4mm,-0.7mm>*{};<4.5mm,-5mm>*{}**@{-},
 <0.6mm,-0.44mm>*{};<8mm,-5mm>*{}**@{-},
 \end{xy}}\Ea
 \Eeq
Here
 \Bi
 \item
 in the first term the differential $\delta$ acts on the vertices of $\Ga$ as in the quotient properad $\weHoQPcd$, that is, by formula (\ref{differential in HoLBcd}) modulo terms involving $(m,n)$-corollas with $m\neq n$.

 \item in the last two terms one attaches  $(m,m)$-corollas, $m\geq 2$,
   to each outgoing leg and each ingoing leg of $\Ga$.
 \Ei

\subsubsection{\bf Remark} The dg prop $\weHoQPcd=\{\weHoQPcd(m,n)\}_{m,n\geq 0}$ contains graphs with no legs at all, i.e.\ its part $\weHoQPcd(0,0)$ does not vanish. Hence one one might consider an extended deformation theory of that prop in which such graphs also ``enter the game", and introduce a further extension, $(\widehat{\HoQP}_{c,d}^{\atop *\circlearrowright})$ of the above prop $(\widehat{\HoQP}_{c,d}^{\atop +\circlearrowright}, \delta^+)$,
 $$
 \widehat{\HoQP}_{c,d}^{\atop +\circlearrowright} \subsetneq \widehat{\HoQP}_{c,d}^{\atop \bigstar \circlearrowright}
 $$
 by adding to the latter a $(0,0)$-generator $\bu$ in degree $1+c+d$ on which the differential acts trivially.
 The extended derivation complex of $\Der(\weHoQPcd)$ can then be defined as follows (cf.\ \cite{AM}),
 $$
 \Der^*(\weHoQPcd) :=  \Der(\widehat{\HoQP}_{c,d}^{\atop \bigstar\circlearrowright}
\rar \widehat{\HoQP}_{c,d}^{\atop \circlearrowright})
 $$
We have obviously a canonical isomorphism of graded vector spaces,
 \Beq\label{2: Der^*(weHoQP) as vector space}
\Der^*(\weHoQPcd)=\prod_{m\geq 0} \left(\widehat{\HoQP}_{c,d}^{\atop \circlearrowright}(m,m) \otimes \sgn_m^{\ot |c|}\otimes \sgn_n^{\ot |d|}\right)^{\bS_m\times \bS_m}[1+c(1-m)+d(1-n)],
\Eeq
i.e.\ the only difference from (\ref{2: Der(weHoQP) as vector space}) is the presence
of the summand with $m=0$. Contrary to the analogous situation studied in \cite{AM} this extension does {\em not}\, give us something really novel as the new extended complex is the {\em direct}\, sum of already introduced complexes,
\Beq\label{2: Der*Hoqois}
\Der^*(\weHoQPcd)=\widehat{\HoQP}_{c,d}^{\atop \circlearrowright}(0,0)[1+c+d] \ \ \oplus \ \ \Der(\weHoQPcd).
\Eeq
Nevertheless it is useful sometimes (see \S 3 below) to consider the extended complex
while studying the most interesting for us deformation complex $\Der(\weHoQPcd)$.

\subsection{Remark}\label{2: remark on star props} Note that representations of the dg prop $\HoLB_{0,1}$ in a vector space $\R^N$
describe Poisson structures (\ref{1: pi formal power series}) which {\em vanish}\, at $0\in \R^N$. This is a restriction which is desirable to avoid, and this can be easily done via introducing a new family of
dg props $\HoLB_{c,d}^{\atop \bigstar}$ as the free prop generated by $(m,n)$-corollas (\ref{2: symmetries of HoLiebcd corollas}) for all possible values $m,n\geq 0$ and equipped with an obvious extension of the differential (\ref{differential in HoLBcd}). Representations of the (wheeled closure)
of the prop $\HoLB_{0,1}^{\atop \bigstar}$ in a dg (finite-dimensional) vector space $V$ are in 1:1 correspondence with arbitrary $\Z$-graded Poisson structures on $V^*$.

\sip

There is an epimorphism of dg props
$$
p:\HoLB_{c,d}^{\atop \bigstar \circlearrowright} \lon \HoQP_{c,d}^{\atop \bigstar \circlearrowright}
$$
so that the above complex $\Der^*(\weHoQPcd)$ can be identified with the deformation complex of that epimorphism (up to a degree shift),
$$
\Der^*(\weHoQPcd)=\Der(\HoLB_{c,d}^{\atop \bigstar \circlearrowright} \rar \HoQP_{c,d}^{\atop \bigstar \circlearrowright})= \Def(\HoLB_{c,d}^{\atop \bigstar \circlearrowright} \rar \HoQP_{c,d}^{\atop \bigstar \circlearrowright})[1].
$$
In particular, the epimorphism $p$ induces a morphism of dg Lie algebras
$$
\Der(\HoLB_{c,d}^{\atop \bigstar \circlearrowright}):=\Der(\HoLB_{c,d}^{\atop \bigstar \circlearrowright}\stackrel{\Id}{\rar} \HoLB_{c,d}^{\atop \bigstar \circlearrowright}) \lon \Der^*(\weHoQPcd).
$$
The derivation complex $\Der(\HoLB_{c,d}^{\atop \bigstar \circlearrowright})$ has been studied in
\cite{AM} where it has been proven that its cohomology is determined (up one rescaling class) by the Kontsevich graph complex $\GC_{c+d+1}^{\geq 2}$ whose definition is recalled below.

\subsection{Directed and undirected graph complexes}
  A {\em graph without legs}, or simply {\em a graph}, $\Ga$ is a 1-dimensional $CW$ complex whose 0-cells are called {\em vertices}\, and 1-cells are called {\em edges}; its set of vertices is denoted by $V(\Ga)$ while the set of  edges by $E(\Ga)$. A graph $\Ga$ is called {\em directed}\, if each edge $e$ comes equipped with a fixed orientation (one of the two possible on an interval). A vertex $v$ of a directed graph is said to have type $(m,n)$ if it has
$m\geq 0$ outgoing edges and $n\geq 0$ incoming edges;
we write in this case $|v|_{in}=n$, $|v|_{out}=m$; the number $|v|:=|v|_{in}+|v|_{out}$ is called the {\em valency}\, of $v$. A $(1,1)$-vertex is called {\em passing}.

\sip

Given any integer $d\in \Z$, we can associate to any graph $\Ga$ with, say, $k$ edges and $l$ vertices,
$$
E(\Ga)=\{e_1, \ldots, e_k\}, \ \ \ \  V(\Ga)=\{v_1,\ldots, v_l\}
$$
a 1-dimensional vector space
$$
\K_\Ga:=\left\{\Ba{ll}  \wedge^{\# E(\Ga)} \K[E(\Ga)] & \text{if}\ d\in 2\Z\\
\wedge^{\# V(\Ga)} \K[V(\Ga)] & \text{if}\ d\in 2\Z+1
\Ea
\right.
$$
where  $\K[E(\Ga)]$ (resp.\, $\K[V(\Ga)]$)   is the linear span of the set of edges (resp.\ vertices); there are at most two different bases (which differ by a sign) of this space given by simple vectors of the form
$$
e_{i_1}\wedge \ldots \wedge e_{i_k}, \ \ \text{respectively},\ \
v_{j_1}\wedge \ldots \wedge v_{j_l}.
$$
An {\em orientation}\, on $\Ga$ is, by definition, a choice of a particular simple basis $or$ of $\K_\Ga$; equivalently, an orientation on $\Ga$ is a choice of ordering of edges (resp.\ vertices)  up to an even permutation. If $\# E(\Ga)\geq 2$
for $d$ even or  $\# V(\Ga)\geq 2$
for $d$ odd, there are precisely two different orientations, $(or, -or)$, on $\Ga$.

\sip

Let $\dfGC_d$ be the ("directed full") completed\footnote{We mean the completion with respect to the filtration of $\dfGC_d$ by the number of vertices.} topological vector space generated over a field $\K$ by the set of all pairs $(\Ga, or)$ (which is often abbreviated simply by $\Ga$) modulo the equivalence relation,
$$
(\Ga, -or) = - (\Ga, or).
$$
We make  $\dfGC_d$ into a $\Z$-{\em graded}\, vector space by setting the (cohomological) degree
of any graph generator $\Ga$ to be given by
$$
|\Ga|= d\# V(\Ga) + (1-d) \# E(\Ga) - d.
$$
 This graded vector space has a Lie algebra structure with
$$
[\Ga',\Ga'']:= \sum_{v\in V(\Ga')} \Ga'\circ_v \Ga'' - (-1)^{|\Ga'||\Ga''|}
\Ga''\circ_v \Ga'
$$
where  $\Ga'\circ_v \Ga''$ is defined by substituting the graph $\Ga''$ into the vertex $v$ of $\Ga'$ and taking the sum over all possible re-attachments of $|v|_{in}+ |v_{out}|$ dangling edges to the vertices of $\Ga''$.
It is easy to see that graph $\xy
 (0,0)*{\bullet}="a",
(5,0)*{\bu}="b",
\ar @{->} "a";"b" <0pt>
\endxy\in \dcGC_d$ is a Maurer-Cartan element, so that one makes $\dfGC_d$ into a complex with the differential
 $$
 \delta:= [\xy
 (0,0)*{\bullet}="a",
(5,0)*{\bu}="b",
\ar @{->} "a";"b" <0pt>
\endxy ,\ ].
 $$
 This dg Lie algebra was introduced and studied in \cite{Ko,Wi1}. The complex $\dfGC_d$ contains a  subcomplex $\dfcGC_d$ spanned by connected graphs. One has an isomorphism of complexes
$$
 \dfGC_d=\odot^{\bu\geq 1} \left(\mathsf{dfcGC}_d [-d]\right)[d]
$$
 The complex $\dfcGC_d$ contains a  subcomplex $\dfcGC_d^{\geq 2}$ spanned by graphs with all vertices $v$ having valency $|v|\geq 2$ which in turn contains a subcomplex $\dcGC_d$
 spanned by graphs with no passing vertices (in particular, with no the tadpole graph consisting of one vertex and one edge); the monomorphism $\dfcGC_d^{\geq 2} \hook \dfcGC_d$ is a quasi-isomorphism, while the monomorphism
 $\dcGC_d\hook \dfcGC_d^{\geq 2}$
 is a quasi-isomorphism up to the tadpole graph \cite{Wi1}.
 The complex $\dcGC_d$  contains a subcomplex $\OGC_d$
  spanned by directed connected graphs with no closed paths of directed edges (``wheels"); such
directed graphs are called {\em oriented}. It was proven in \cite{Wi2} that
$$
H^\bu(\dcGC_d) =H^\bu(\OGC_{d+1}).
$$

\sip

Let $\dfGC_d^{(l)} \subset \dfGC_d$ be the linear subspace spanned by graphs with $l$ edges.
The group $\Z_2^{\times l}$ acts naturally on the generators of $\dfGC_d^{(l)}$   by reversing
directions on edges. Hence one can define an undirected version of $\dfGC_d$ as follows,
$$
\mathsf{fGC}_d= \prod_{l\geq 0}\dfGC_d^{(l)} \ot_{ \Z_2^{\times l} } \sgn_2^{\ot l|d|}
$$
Graphs from $\fGC_d$ have no directions on edges for even $d$, or have a direction  given on every edges up a flip and multiplying by $-1$.

\sip

The complex $\mathsf{fGC}_d$ contains a subcomplex $\fGC_d^{\geq 2}$ spanned by graphs with all vertices having valency at least two, which in turn contains a subcomplex $\GC_d^{\geq 2}$ spanned by connected graphs. There is an isomorphism of complexes,
\Beq\label{2: fGC versus GC}
 \fGC_d^{\geq 2}=\odot^{\bu\geq 1} \left(\mathsf{GC}^{\geq 2}_d [-d]\right)[d]
 \Eeq
and
a direct splitting,
$$
\GC_d^{\geq 2}= \GC_d \ \oplus\ \GC_d^{2}
$$
where $\GC_d^{2}$ is the subcomplex spanned by graphs containing at least one bivalent vertex and  $\GC_d$ is spanned by graphs with all vertices at least trivalent. The cohomology of $\GC_d^{2}$ has been computed  in \cite{Wi1},
$$
H^\bu(\GC_d^2)=
\bigoplus_{p\geq 1\atop p\equiv 2d+1 \mod 4} \K[d-p],
$$
where the summand $\K[d-p]$ is generated by the polytope with $p$ vertices,
that is, a connected graph with $p$ bivalent vertices. On the other hand, the cohomology of $\GC_d$ is understood at present only in non-positive degrees \cite{Wi1}
$$
H^{\bu \leq -1}(\GC_2)=0\ \ \text{and}\ \ H^0(\GC_2)=H^0(\GC_2^{\geq 2})=H^0(\OGC_3)=\grt_1,
$$
where
$\grt_1$ is the Lie algebra of the Grothendieck-Teichm\"uller group $GRT_1$. Moreover, it was proven in \cite{Wi1} that there is a quasi-isomorphism of complexes
\Beq\label{2: from GC to dcGC}
i: \GC_d^{\geq 2} \lon \dfcGC_d^{\geq 2}
\Eeq
which sends an undirected graph into a sum of directed graphs by assigning directions to its edges
in all possible ways.

\sip

In this paper we have to work with the ``full" Grothendieck-Teichm\"uller group $GRT$ rather than with its subgroup $GRT_1$. Hence we need a slight extension \cite{MW2} of the complex $\GC_d^{\geq 2}$,
$$
\GC_{d}^{\geq 2} \lon \GC_{d}^{\geq 2}\oplus \K
$$
where the generator of the 1-dimensional summand $\K$  is a cycle, ``a graph with no vertices and edges", and has Lie bracket with any graph $\Ga$ from $\GC_{d}^{\geq 2}$ equal to $2\ell \Ga$, where
$\ell$ is the number of its loops. In our story this class takes care about the obvious rescaling automorphism of the prop $\weHoQPcd$. Hence we define the {\em full graph complexe}\, of not necessarily connected graphs as the completed graded symmetric tensor algebra
\Beq\label{1: fGC in terms of GC}
 \mathsf{FGC}^{\geq 2}_{d}:=\widehat{\odot^\bu}\left((\GC_{d}^{\geq 2}\oplus \K )[-d]\right)[d].
\Eeq
%
%
%
The above results imply
$$
H^0(\mathsf{FGC}_2^{\geq 2})=\grt,
$$
where $\grt$ is the Lie algebra of $GRT$.

\sip

\subsection{A quotient complex of $\dcGC_d$}\label{2: section on dcGC=}
The complex $\dcGC_d$ contains a subcomplex $\dcGC_{d}^{\neq}$ spanned by graphs with at least one vertex $v$ satisfying $|v|_{in}\neq |v|_{out}$. The quotient complex $\dcGC_d^{=}$,
\Beq\label{2: short ex seq for H(dGC=)}
0\lon \dcGC_{d}^{\neq} \lon \dcGC_d \lon \dcGC_d^{=} \lon 0,
\Eeq
is spanned by graphs with all vertices satisfying $|v|_{in}=|v|_{out}$.
It is precisely the connected part of the direct summand $\widehat{\HoQP}_{c',d'}^{\atop  \circlearrowright}(0,0)[1+c'+d']$ in the ``full" derivation complex (\ref{2: Der*Hoqois}) with $d=c'+d'+1$.



\sip

 Any graph $\Ga$ in $\dcGC_d^{=}$ has vertices $v\in V(\Ga)$ of type $(l_v,l_v)$ for some $l_v\geq 2$; hence its cohomological degree is given by
$$
|\Ga|=d(\# V(\Ga)-1) + (1-d)\sum_{v\in V(\Ga)} l_v.
$$
Hence for $d\geq 2$ we obtain
$$
|\Ga| \leq d(\# V(\Ga)-1) + (1-d)\sum_{v\in V(\Ga)}2= (2-d)\# V(\Ga)-d
$$
and conclude that for $d=2$ the quotient complex $\dcGC_2^{=}$ is concentrated in degrees $\leq -2$. As $H^\bu(\dcGC_2)$ is concentrated in non-negative degrees \cite{Wi1}, we also conclude that the canonical projection at the cohomology level
$$
H^\bu(\dcGC_2) \lon H^\bu(\dcGC_2^{=})
$$
 is equal to the zero map (hence the map $H^\bu(\dcGC_d) \lon H^\bu(\dcGC_d^{=})$  is also equal to zero for any {\em even}\, $d$). Therefore in the most important for us in this paper  case $d=2$ the long exact sequence of cohomology groups associated with
(\ref{2: short ex seq for H(dGC=)}) decomposes into a collection of short exact sequences,
\Beq\label{2: short ex seq for H(dGC=2)}
0 \lon H^{\bu}(\dcGC_2^{=}) \stackrel{b}{\lon} H^{\bu +1}(\dcGC_2^{\neq}) \lon H^{\bu+1}(\dcGC_2) \lon 0
\Eeq
The subcomplex $\widehat{\HoQP}_{0,1}^{\atop \circlearrowright}(0,0)[2]=\odot^{\geq 1}(\dcGC_2^{=}[-2])[2]$  in the ``full" derivation complex (\ref{2: Der*Hoqois}) does not interact with the second summand in
(\ref{2: Der*Hoqois}) at all; it  is the second summand $\Der(\weHoQPcd)$
which plays the key role  in classification of the universal quantizations of $\Z$-graded Poisson structures. We shall study its cohomology in \S 4 below.

\bip

\bip


{\Large
\section{\bf Cohomology of the deformation complex of $\HoQP_{c,d}^\circlearrowright$}
}

\sip

\subsection{\bf Graph complexes versus derivation complexes}
The differentials in the dg props $\HoLBcd$, $\HoLBcd^{\atop \bigstar}$ and $\HoQP_{c,d}$ (and of their wheeled closures)
are given by connected graphs (see (\ref{differential in HoLBcd}) so the subspaces of all these props spanned by {\em connected}\, graphs form {\em dg}\ properads which we denote by the same symbols. However we use a slightly different symbol $\sDer$ to denote the Lie algebras of derivations of the properads as opposed to the symbol $\Der$ is reserved for derivations of the associated props. Both types of dg Lie algebras are related to each other by an exact symmetric tensor algebra functor as follows,
\Beqrn
\Der(\HoLBcd)&=&\odot^{\bu\geq 1}(\sDer(\HoLBcd)[-1-c-d])[1+c+d], \\ \Der(\wHoQP_{c,d})&=&\odot^{\bu\geq 1}(\sDer(\wHoQP_{c,d})[-1-c-d])[1+c+d]
\Eeqrn
Thus it is enough to study the cohomology of the complex $\sDer$  of {\em properadic}\, derivations in each case.

\sip

There is a morphism of dg Lie algebras (see \S 3.3 in \cite{MW1})
\Beq\label{2: Morhism F from GC_3^or}
\Ba{rccc}
 F\colon & \OGC_{c+d+1}  &\to & \sDer(\HoLBcd)\\
         &   \Ga & \to & F(\Ga)
         \Ea
\Eeq
where the derivation $F(\Ga)$ acts (from the right) on the generators of the  completed (by the loop number, cf. \cite{MW1}) properad  $\wHoLBcd$ as follows
\Beq \label{equ:def GC action 1}
\left(\Ba{c}\resizebox{12mm}{!}{\begin{xy}
 <0mm,0mm>*{\circ};<0mm,0mm>*{}**@{},
 <-0.6mm,0.44mm>*{};<-8mm,5mm>*{}**@{-},
 <-0.4mm,0.7mm>*{};<-4.5mm,5mm>*{}**@{-},
 <0mm,0mm>*{};<-1mm,5mm>*{\ldots}**@{},
 <0.4mm,0.7mm>*{};<4.5mm,5mm>*{}**@{-},
 <0.6mm,0.44mm>*{};<8mm,5mm>*{}**@{-},
   <0mm,0mm>*{};<-8.5mm,5.5mm>*{^1}**@{},
   <0mm,0mm>*{};<-5mm,5.5mm>*{^2}**@{},
   <0mm,0mm>*{};<9.0mm,5.5mm>*{^m}**@{},
 <-0.6mm,-0.44mm>*{};<-8mm,-5mm>*{}**@{-},
 <-0.4mm,-0.7mm>*{};<-4.5mm,-5mm>*{}**@{-},
 <0mm,0mm>*{};<-1mm,-5mm>*{\ldots}**@{},
 <0.4mm,-0.7mm>*{};<4.5mm,-5mm>*{}**@{-},
 <0.6mm,-0.44mm>*{};<8mm,-5mm>*{}**@{-},
   <0mm,0mm>*{};<-8.5mm,-6.9mm>*{^1}**@{},
   <0mm,0mm>*{};<-5mm,-6.9mm>*{^2}**@{},
   <0mm,0mm>*{};<9.0mm,-6.9mm>*{^n}**@{},
 \end{xy}}\Ea\right)\cdot F(\Ga)
=
 \sum_{s:[n]\rar V(\Ga)\atop \hat{s}:[m]\rar V(\Ga)}  \Ba{c}\resizebox{9mm}{!}  {\xy
 (-6,7)*{^1},
(-3,7)*{^2},
(2.5,7)*{},
(7,7)*{^m},
(-3,-8)*{_2},
(3,-6)*{},
(7,-8)*{_n},
(-6,-8)*{_1},
(0,4.5)*+{...},
(0,-4.5)*+{...},
(0,0)*+{\Ga}="o",
(-6,6)*{}="1",
(-3,6)*{}="2",
(3,6)*{}="3",
(6,6)*{}="4",
(-3,-6)*{}="5",
(3,-6)*{}="6",
(6,-6)*{}="7",
(-6,-6)*{}="8",
\ar @{-} "o";"1" <0pt>
\ar @{-} "o";"2" <0pt>
\ar @{-} "o";"3" <0pt>
\ar @{-} "o";"4" <0pt>
\ar @{-} "o";"5" <0pt>
\ar @{-} "o";"6" <0pt>
\ar @{-} "o";"7" <0pt>
\ar @{-} "o";"8" <0pt>
\endxy}\Ea
\Eeq
 where the sum is taken over all ways of attaching the incoming and outgoing legs to the graph $\Ga$, and one sets to zero every resulting graph if it contains (i) a vertex with valency $<3$, or (ii)  or a vertex with no incoming edge(s)/leg(s), or (iii) a vertex with no outgoing edge(s)/leg(s).
 Moreover, it was proven in \cite{MW1} that map $F$ is a quasi-isomorphism up to one class in the complex $\sDer(\HoLBcd)$ given explicitly by the cycle,
 \Beq\label{4: rescaling class r}
 r=\displaystyle
  \sum_{m,n}(n-1)
  \overbrace{
  \underbrace{
 \Ba{c}\resizebox{6mm}{!}  {\xy
(0,4.5)*+{...},
(0,-4.5)*+{...},
(0,0)*{\bu}="o",
(-5,5)*{}="1",
(-3,5)*{}="2",
(3,5)*{}="3",
(5,5)*{}="4",
(-3,-5)*{}="5",
(3,-5)*{}="6",
(5,-5)*{}="7",
(-5,-5)*{}="8",
\ar @{-} "o";"1" <0pt>
\ar @{-} "o";"2" <0pt>
\ar @{-} "o";"3" <0pt>
\ar @{-} "o";"4" <0pt>
\ar @{-} "o";"5" <0pt>
\ar @{-} "o";"6" <0pt>
\ar @{-} "o";"7" <0pt>
\ar @{-} "o";"8" <0pt>
\endxy}\Ea
 }_{n\times}
 }^{m\times}.
 \Eeq
Note that coefficient $(n-1)$ above can be replaced by $(m-1)$  or by $(m+n-2)$ --- all the corresponding cycles represent the same cohomology class.

 \sip

A  similar morphism of dg Lie algebras
$$
F^\circlearrowright: \dcGC_{c+d+1} \lon \sDer({\HoLB}_{c,d}^{\atop \bigstar\circlearrowright})
$$
for the wheeled closure  ${\HoLB}_{c,d}^{\atop \bigstar\circlearrowright}$ of the ``full" ordinary dg properad
has been studied in \cite{AM}
where it has been proven that $F^\circlearrowright$ is also a quasi-isomorphism up to the same rescaling class $r$; that proof is in fact much easier and shorter than the proof of the quasi-isomorphism (\ref{2: Morhism F from GC_3^or})) in the unwheeled case.

\sip

The first (resp. second) result has been used in \cite{MW3} (resp. in \cite{AM}) to classify (up to homotopy) all universal deformation quantizations of Lie bialgebras (resp., of generic Poisson structures).

\sip

In this section we address a similar problem for the derivation complex of the wheeled properad of homotopy quadratic Poisson structures. Its solution given below will lead us in \S 5 to the classification of all homotopy classes of universal quantizations of $\Z$-graded quadratic Poisson structures.

\subsection{Marked graph complex}
A {\em weight function}\, on a directed graph $\Ga\in \dfcGC_d$
is, by definition,  a map
$$
\Ba{rccc}
\fw: & V(\Ga) &\lon &  \N_{\geq 1}\\
              & v & \lon & \fw_v
\Ea
$$
satisfying the condition
    \Beq\label{3: condition on weights}
    \fw_v\geq \max\{1,|v|_{in}-1, |v|_{out}-1\}
    \Eeq
The weight function $\fs$ is given by the equality
$$
\fs_v:=\max\{1,|v|_{in}-1, |v|_{out}-1\}
$$
is called {\em canonical}. 

\sip

A {\em weighted}\, graph is, by definition, a pair $\Ga^\fw:=(\Ga, \fw)$ where $\Ga\in \dfcGC_d$ and $\fw$ is a weight function on $\Ga$.
The integer $\fw_v$ is called the {\em weight}\, of the vertex $v$ of a weighted graph. The natural number
$
\fw(\Ga):= \sum_{v\in V(\Ga)} \fw_v
$
is called the {\em (total) weight}\, of the weighted graph $\Ga^\fw$. The cohomological degree of $\Ga^\fw$ is defined by the standard formula
\Beq\label{4: |Ga^w|}
|\Ga^\fw|:= d(\# V(\Ga)-1) - (d-1) \# E(\Ga)
\Eeq
i.e.\ $|\Ga^\fw|=|\Ga|$, the degree of the underlying graph.


\sip

Let $\WGC_d^\star$ be the completed (by the number of vertices and edges) graded vector space over a field $\K$ which is generated by
weighted directed connected graphs. Its generators $\Ga^\fw$ can be represented pictorially as follows,
$$
\Ba{c}
\resizebox{15mm}{!}{ \xy
(0,0)*+{_a}*\cir{}="d1",
(10,0)*+{_b}*\cir{}="d2",
(-5,-5)*{}="dl",
(5,-5)*{}="dc",
(15,-5)*{}="dr",
(0,10)*+{_c}*\cir{}="u1",
(10,10)*+{_d}*\cir{}="u2",
(5,15)*{}="uc",
(5,15)*{}="uc",
(15,15)*{}="ur",
(0,15)*{}="ul",
\ar @{<-} "d1";"d2" <0pt>
\ar @{<-} "u1";"d1" <0pt>
\ar @{->} "u1";"u2" <0pt>
\ar @{<-} "u1";"d2" <0pt>
\ar @{->} "u2";"d2" <0pt>
\ar @{<-} "u2";"d1" <0pt>
\endxy}
\Ea
,
\ \ \ \
\Ba{c}
\resizebox{15mm}{!}{ \xy
(0,0)*+{_a}*\cir{}="d1",
(10,0)*+{_b}*\cir{}="d2",
(-5,-5)*{}="dl",
(5,-5)*{}="dc",
(15,-5)*{}="dr",
(0,10)*+{_c}*\cir{}="u1",
(10,10)*+{_d}*\cir{}="u2",
(5,15)*{}="uc",
(5,15)*{}="uc",
(15,15)*{}="ur",
(0,20)*+{_e}*\cir{}="ul",
\ar @{<-} "d1";"d2" <0pt>
\ar @{<-} "u1";"d1" <0pt>
\ar @{->} "u1";"u2" <0pt>
\ar @{<-} "u1";"d2" <0pt>
\ar @{->} "u2";"d2" <0pt>
\ar @{<-} "u2";"d1" <0pt>
\ar @{<-} "ul";"u1" <0pt>
\endxy}
\Ea
 \in \WGC_d^\star,
$$
where the natural numbers $a,b,d,c,e\in \N_{\geq 1}$ stand for weights $\fw_v$ of the corresponding vertices.

\sip
Recall that the standard differential $\delta$ in $\dfcGC_d$ is defined by
\Beq\label{4: delta in dfcGC}
\delta \Ga:= [\Ga, \ \xy
 (0,-2)*{\bullet}="a",
(0,2)*{\bu}="b",
\ar @{->} "a";"b" <0pt>
\endxy\ ]=
\delta_0\Ga - (-1)^{|\Ga|}  \sum_{v\in V(\Ga)}
\Ba{c}
{ \xy
(0,7)*{\bu}="u1",
(0,0)*+{\Ga}="ul",
\ar @{->} "ul";"u1" <0pt>
\endxy}
\Ea
- (-1)^{d+|\Ga|} \sum_{ v\in V(\Ga)}
\Ba{c}
{ \xy
(0,0)*{\bu}="u1",
(0,7)*+{\Ga}="ul",
\ar @{<-} "ul";"u1" <0pt>
\endxy}
\Ea
\Eeq
where
\Bi
\item[-] the first summand $\delta_0 \Ga$ is equal to $\sum_{v\in V(\Ga)}\delta_v$ with $\delta_v\Ga$ obtained
from $\Ga$ by replacing
   the vertex $v$ with the graph   $\xy
 (0,-2.6)*{\bullet}="a",
(0,2.6)*{\bu}="b",
\ar @{->} "a";"b" <0pt>
\endxy$
   and taking the sum over all possible ways to reattach  the dangling  edges (connected earlier to $v$) to the two newly created vertices.

\item[-] the second (resp., the third) summand is obtained by attaching the ingoing (resp., outgoing) leg of \ $\xy
 (-2,0)*{\bullet}="a",
(2,0)*{\bu}="b",
\ar @{->} "a";"b" <0pt>
\endxy$ to a vertex $v\in V(\Ga)$.
\Ei

The graphs $\Ga$ in $\dfGC_d$ are oriented (depending on the parity of $d$), and their orientations are defined up to a sign. Let us explain in full detail the rule of signs behind formula (\ref{4: delta in dfcGC}) by defining explicitly the induced orientation of each summand on the r.h.s.\ of (\ref{4: delta in dfcGC}) from the given orientation $or$ of an input graph $\Ga$:

\Bi
\item[(i)] if $d\in 2\Z$, then $or$ is given by choosing of an ordering of edges, say, by
$$
or=e_1\wedge e_2\wedge \ldots \wedge e_k, \ \ \ \ k=\# E(\Ga).
$$
Then each graph in the first summand $\delta_{0}\Ga$ has, by definition, the orientation given by
$$
e_1\wedge e_2\wedge \ldots \wedge e_k\wedge e'
$$
where $e'$ is the newly created edge. Each graph in the remaining two summands (the ones with the univalent black vertex) has an induced orientation given by
$$
e'\wedge e_1\wedge e_2\wedge \ldots \wedge e_k.
$$
As $(-1)^{|\Ga|}=(-1)^k$, we conclude that if $\Ga$ is not an isolated vertex, then the two  graphs in the summand $\delta_{v_i}\Ga$ of $\delta_0\Ga$  which have a univalent black vertex cancel out the corresponding graphs from the remaining two summands which have the new edge attached to $v_i$.

\item[(ii)] if $d\in 2\Z+1$, then $or$ is given by choosing an ordering of vertices, say
$$
or=v_1\wedge v_2\wedge\ldots \wedge v_i\wedge  \ldots \wedge v_l, \ \ \ \ l=\# V(\Ga).
$$
Then each graph in the summand $\delta_{v_i}\Ga$ of $\delta_0\Ga$  has, by definition, the orientation given by
$$
 v_1\wedge \ldots \wedge v_{i-1}\wedge v_i'\wedge v_{i+1}\wedge \ldots \wedge v_l\wedge v_i'',
$$
where $v_i'$ and $v_i''$ are the vertices of the newly created edge $e'$ such that $v_i'$ (resp., $v_i''$) is the initial (resp.\ last) vertex along the direction flow on $e'$. Each  graph in the last two summands in (\ref{4: delta in dfcGC}) which has the new edge attached to $v_i$ (in one of the two possible ways) has the induced orientation given by
$$
v_i'\wedge v_1\wedge \ldots\wedge  v_{i-1}\wedge v_i''\wedge v_{i+1}\wedge \ldots \wedge v_l
$$
As $(-1)^{|\Ga|}=(-1)^{l+1}$, we conclude that if $\Ga$ is not an isolated vertex, then the two  graphs in each summand $\delta_{v_i}\Ga$ which have $v'$ or $v''$ univalent cancel out the corresponding graphs from the remaining two summands, the ones have the new edge attached to $v_i$.

\Ei

\sip
We make $\WGC_d$ into a complex  by defining a differential $d$ in a close analogy to the above formula,
 \Beq\label{4: d in wcGC}
 d\Ga^\fw =d_0 \Ga^\fw - (-1)^{|\Ga|} \sum_{c\in \N_{\geq 1}}\sum_{ v\in V(\Ga)\atop
 \fw_v> |v|_{out}-1}
\Ba{c}
{ \xy
(0,8)*+{_c}*\cir{}="u1",
(0,0)*+{\Ga^\fw}="ul",
\ar @{->} "ul";"u1" <0pt>
\endxy}
\Ea
-  (-1)^{|\Ga|+d} \sum_{c\in \N_{\geq 1}}\sum_{ v\in V(\Ga)\atop
 \fw_v> |v|_{in}-1}
\Ba{c}
{ \xy
(0,0)*+{_c}*\cir{}="u1",
(0,8)*+{\Ga^\fw}="ul",
\ar @{<-} "ul";"u1" <0pt>
\endxy}
\Ea
 \Eeq
where
\Bi
\item[-] the first summand $d_0 \Ga^\fw$ is equal to $\sum_{v\in V(\Ga)}\delta_v\Ga^\fw$ with $\delta_v\Ga^\fw$ obtained
from $\Ga^\fw$ by replacing
   the vertex $v$ with weight, say, $a\in \N_{\geq 1}$ by the sum of graphs
   \Beq\label{5: splitting weight vertex}
 \xy
(0,0)*+{_{a}}*\cir{}="u1",
\endxy \ \ \lon \ \ \ \   \sum_{a=a'+ a''}
\Ba{c}   \resizebox{4mm}{!}
   {\xy
(0,10)*+{_{a'}}*\cir{}="u1",
(0,0)*+{_{a''}}*\cir{}="ul",
\ar @{->} "ul";"u1" <0pt>
\endxy}\Ea
   \Eeq
   and taking one more sum over all possible way to reattach  the dangling  edges to the newly created vertices $\Ba{c}   \resizebox{4mm}{!}
   {\xy
(0,1)*+{_{a'}}*\cir{},
\endxy}\Ea$ and $\Ba{c}   \resizebox{4mm}{!}
   {\xy
(0,1)*+{_{a''}}*\cir{},
\endxy}\Ea$ in such a way
   that the condition (\ref{3: condition on weights}) is respected; for example, if $v$ has weight $1$, then $\delta_v \Ga^\fw=0$;

\item[-] the second summand consists of graphs obtained by attaching the ingoing leg of the
graph $\Ba{c}\resizebox{3mm}{!}{\xy
(0,6)*+{_c}*\cir{}="u1",
(0,0)*{}="ul",
\ar @{->} "ul";"u1" <0pt>
\endxy}\Ea$ to each vertex $v$ of $\Ga^\fw$ whose weight satisfies the strict inequality
$\fw_v>\# |v|_{out}-1$ (so that the condition (\ref{3: condition on weights}) is not violated);

\item[-] the third summand is an obvious ``outgoing" analogue of the second one.
\Ei
The rule of signs in (\ref{4: d in wcGC}) is identical to the rule of signs in (\ref{4: delta in dfcGC}) as weights on vertices do not influence the cohomological degrees and orientations of weighted graphs! This simplification of the rule of signs is one of the motivation of introducing the complex $\WGC_d$.

\sip

This complex splits as a direct sum (cf.\ (\ref{2: Der*Hoqois}))
\Beq\label{5: Decomp od WGC^star}
\WGC_d^\star=\WGC_d^{=} \oplus \WGC_d
\Eeq
where the summand $\WGC_d^{=}$ is spanned by graphs $\Ga$ such that for each $v\in V(\Ga)$ one has $\fw_v=|v|_{in}-1=|v|_{out}-1$ (implying $|v|_{in}=|v|_{out}\geq 2$); the differential in this summand is given by the first term $d_0$ only in the above general formula (\ref{4: d in wcGC})  for $d$. It is not hard to check that
$\WGC_d^{=}$ is identical to the quotient complex $\dcGC_d^=$ of $\dcGC_d$ introduced in \S
{\ref{2: section on dcGC=}}. This observation tells us that the similarity between formulae
(\ref{4: delta in dfcGC}) and (\ref{4: d in wcGC}) can be sometimes deceiving --- the complex
$\dcGC_d^=$ is {\em not}\, a direct summand of $\dcGC_d$.

\sip

The summand $\WGC_d$
is spanned by graphs with at least one vertex $v$ having weight $\fw_v\neq |v|_{in}-1$ or $\fw_v\neq |v|_{out}-1$.

\subsubsection{\bf Lemma} {\em There are  canonical isomorphisms of complexes,
$$
\sDer^\star(\wHoQPcd)\equiv \sDer(\HoLB_{c,d}^{\atop \bigstar \circlearrowright} \rar \HoLB_{c,d}^{\atop \circlearrowright})    \simeq  \WGC_{c+d+1}^\star, \ \ \ \ \  \sDer(\wHoQPcd)\simeq  \WGC_{c+d+1}.
$$
}
\begin{proof}
Consider a map
$$
\Ba{rccc}
w: & \sDer^\star(\wHoQPcd) &\lon & \WGC_{c+d+1}^\star\\
& \ga &\lon &\ga^w
\Ea
$$
 defined as follows: $\ga^w$ is obtained from $\ga\in \sDer^*(\wHoQPcd)$ by decorating each $(m,m)$-vertex $v$ of $\ga$, $m\geq 2$, with the weight $\fw_v:=m-1$ and removing all ingoing and outgoing legs. This map is an epimorphism: given any $\ga^w$, the graph $\ga$ obtained from it by attaching to each vertex $v\in V(\ga^w)$ $\fw_v+1 - |v_{out}|$ outgoing legs and $\fw_v+1 - |v_{in}|$ ingoing ones (and then (skew)symmetrizing them) belongs to the pre-image $w^{-1}(\ga^w)$; moreover, this $\ga$ is reconstructed from
$\ga^w$ is reconstructed from $\ga^w$ in a unique way, so the map $w$ is an isomorphism
of vector spaces.

\sip

Let us check next that the map $w$ respects the $\Z$-grading of both sides. On the one hand we have (cf.\ (\ref{4: |Ga^w|})
$$
|\ga^w|= (c+d+1)(\# V(\ga)-1) - (c+d)\# E(\ga).
$$
On the other hand, the isomorphism (\ref{2: Der(weHoQP) as vector space}) implies that each vertex $v\in V(\ga)$,
that is, each $(\fw_v+1,\fw_v+1)$-corolla in $\ga$ contributes the degree
$$
1+ c+d - (c+d)(a_v +1)= 1-(c+d)\fw_v
$$
Shifting the resulting total degree by $(1+c+d -(c+d)\frac{1}{2} L(\ga))$, where $L(\ga)$
is the total number of ingoing and outgoing legs of $\ga$, we obtain
\Beqrn
|\ga| &=& \sum_{v\in V(\ga)} \left( 1-(c+d)\fw_v\right) - \left(1+c+d -(c+d)\frac{1}{2} L(\ga)\right)\\
      &=& \# V(\Ga) - (c+d)\fw(\ga^w) - \left(1+c+d -(c+d)\frac{1}{2} L(\ga)\right).
\Eeqrn
The total number of half-edges in $\Ga$ is equal to $2\sum_{v\in V(\ga)} (1+\fw_v)=2\# V(\ga) +2\fw(\ga^w)$, out of this number $2\# E(\ga)$ are paired into internal edges, so that the remaining number of legs is given by
$$
L(\ga)=2\# V(\ga) +2\fw(\ga) -2\# E(\ga)
$$
and we obtain finally,
\Beqrn
|\ga|&=& \# V(\ga) - (c+d)\fw(\ga^w) - \left(1+c+d -(c+d)\left(\# V(\ga) +\fw(\ga^w) -\# E(\ga)\right)\right)\\
&=& (c+d+1)(\# V(\ga)-1) - (c+d)\# E(\ga).
\Eeqrn
Hence $|\ga^w|=|\ga|$ so that the map $w$ does respect the degrees.

\sip

The compatibility of the map $w$ with the differentials follows from the very definition (\ref{4: d in wcGC}) of the differential in $\WGC_{c+d+1}$ ---  the latter was essentially copied from the differential (\ref{d in Der(Holieb)}) in $\Der^\star(\wHoQPcd)$ via the isomorphism $w$.
\end{proof}

It is easy to check that
\Beq\label{4: weighted rescaling class r}
r:=\sum_{a\geq 1} (a-1)\ \xy
(0,0)*+{_{a}}*\cir{}="u1",
\endxy
\Eeq
is a non-trivial cohomology class in $\WGC_d$. It is the image of the rescaling class (\ref{4: rescaling class r}) under the composition of the projection $\sDer(\wHoLBcd)\rar \sDer(\wHoQPcd)$ with the above isomorphism.

\sip

There is a monomorphism of complexes
\Beq\label{4: F from dcGC to wGC}
\Ba{rccc}
\cF^\circlearrowright:  & \dfcGC^{\geq 2}_d & \lon & \WGC_d\\
   &  \Ga   & \lon & \displaystyle \sum_{\fw: V(\Ga)\rar \N_{\geq 1}} \Ga^{\fw}
\Ea
\Eeq
where $\Ga^\fw$ is the graph obtained from $\Ga$ by assigning to each vertex $v$ the weight  $a_v:=\fw(v)$, and setting $\Ga^\fw$ to zero if for at least one vertex the condition (\ref{3: condition on weights}) is violated.

\sip

\subsection{\bf An auxiliary complex}\label{2: subsec on C}  Consider an auxiliary graph complex,
\Beq\label{aux complex C}
C =\bigoplus_{n\geq 1} C_n.
\Eeq
where  $C_n$ is spanned by graphs of the form
\Beq\label{5: element of C}
\Ba{c}
\resizebox{35mm}{!}
{ \xy
(-20,1)*+{_{a_1}}*\cir{}="1",
(-10,1)*+{_{a_2}}*\cir{}="2",
(-2,1)*{}="3",
(1,1)*{\ldots},
(4,1)*{}="4",
(12,1)*+{_{a_n}}*\cir{}="5",
(20,1)*{}="6",
%
\ar @{->} "1";"2" <0pt>
\ar @{->} "2";"3" <0pt>
\ar @{->} "4";"5" <0pt>
\ar @{->} "5";"6" <0pt>
\endxy}
\Ea,
\Eeq
 with $a_1,\ldots, a_n\in \N_{\geq 1}$. It is useful to interpret such graphs as elements of a free properad generated by weighted corollas of type $(1,1)$ and $(1,0)$ because the differential can defined by its action on generators as follows
\Beq\label{d on C}
d\Ba{c}
\resizebox{8mm}{!}
{ \xy
(12,1)*+{_{a}}*\cir{}="5",
(18,1)*{}="6",
%
\ar @{->} "5";"6" <0pt>
\endxy}
\Ea
=
\sum_{a=b+c}
\Ba{c}
\resizebox{16mm}{!}
{ \xy
(-20,1)*+{_{b}}*\cir{}="1",
(-11,1)*+{_{c}}*\cir{}="2",
(-4,1)*{}="3",
%
\ar @{->} "1";"2" <0pt>
\ar @{->} "2";"3" <0pt>
\endxy}
\Ea,
\ \ \ \
d\Ba{c}
\resizebox{12mm}{!}
{ \xy
(6,1)*{}="4",
(12,1)*+{_{a}}*\cir{}="5",
(18,1)*{}="6",
\ar @{->} "4";"5" <0pt>
\ar @{->} "5";"6" <0pt>
\endxy}
\Ea
=
\sum_{a=b+c}
\Ba{c}
\resizebox{22mm}{!}
{ \xy
(-27,1)*{}="0",
(-20,1)*+{_{b}}*\cir{}="1",
(-11,1)*+{_{c}}*\cir{}="2",
(-4,1)*{}="3",
\ar @{->} "0";"1" <0pt>
\ar @{->} "1";"2" <0pt>
\ar @{->} "2";"3" <0pt>
\endxy}
\Ea,
\Eeq
\subsubsection{\bf Lemma \cite{CMW}}\label{2: lemma on C}
$H^\bu(C)= \K[-1]\simeq \mbox{span} \langle
\Ba{c}
\resizebox{8mm}{!}
{ \xy
(12,1)*+{_{1}}*\cir{}="5",
(18,1)*{}="6",
\ar @{->} "5";"6" <0pt>
\endxy}
\Ea  \rangle$.

\sip

As the cohomology of $C$ is one-dimensional and is concentrated in degree $1$, we notice for the future use that
$H^\bu(\odot^{k\geq 2} C)=0$.

\subsection{Theorem}\label{4: Theorem on dcGC and dwGC} {\em The morphism
 \Beq\label{4: Theorem on F from dGC to wGC}
\Ba{rccc}
\cF^\circlearrowright:  & \dfcGC_d^{\geq 2}  & \lon & \WGC_d\\
   &  \Ga   & \lon & \displaystyle \sum_{\fw: V(\Ga)\rar \N_{\geq 1}} \Ga^{\fw}
\Ea
\Eeq
is a quasi-isomorphism up to one rescaling class (\ref{4: weighted rescaling class r}).}

\sip



\begin{proof} We shall study the full graph complex $\WGC_d^\star$ and, using certain spectral sequences, show the quasi-isomorphism of complexes
$$
\WGC_d^\star\equiv \WGC_d^{=} \oplus \WGC_d \simeq \dcGC^{=} \oplus \cF^\circlearrowright(\dfcGC_d^{\geq 2}) \oplus \text{span}\langle r\rangle.
$$
which implies the theorem.

\sip

Let us call (weighted or un-weighted) univalent vertices of our graphs with precisely {\em one}\,  outgoing (resp., ingoing)  edge {\em univalent in-vertices}\, (resp. {\em univalent out-vertices}); such vertices have the form $\Ba{c}\resizebox{7mm}{!}{\xy
(0,1)*+{_a}*\cir{}="0",
(6,1)*{}="r",
\ar @{->} "0";"r" <0pt>
\endxy}\Ea$ or $\Ba{c}\resizebox{5mm}{!}{\xy
(0,1)*{\bu}="0",
(6,1)*{}="r",
\ar @{->} "0";"r" <0pt>
\endxy}\Ea$
 (resp.  $\Ba{c}\resizebox{7mm}{!}{\xy
(0,1)*+{_a}*\cir{}="0",
(6,1)*{}="r",
\ar @{<-} "0";"r" <0pt>
\endxy}\Ea$ or $\Ba{c}\resizebox{5mm}{!}{\xy
(0,1)*{\bu}="0",
(6,1)*{}="r",
\ar @{<-} "0";"r" <0pt>
\endxy}\Ea$).
\sip

\sip

{\sf Step 1: skeleton filtration}.   Delete from a graph $\Ga^\fw\in \WGC_d^\star$ recursively all univalent in-vertices. The result is a weighted graph $\Ga^\fw_{sk}$ called the {\em skeleton graph}\, of $\Ga^\fw$; the vertices of $\Ga^\fw_{sk}$ are called the {\em skeleton vertices}\, of $\Ga^\fw$. If we forget the induced weight function on $\Ga_{sk}^\fw$ we get a graph $\Ga_{sk}$ which belongs to a subspace $\dfcGC_d'\subset \dfcGC_d$
of the full directed graph complex spanned by graphs with no univalent in-vertices. For a vertex $r\in \Ga_{sk}$ set
$$
\fs_r^{sk}:=\max\{1,  |r|^{sk}_{in} -1,  |r|^{sk}_{out} -1\},
$$
where $ |r|^{sk}_{in}$ (resp., $ |r|^{sk}_{out}$) is the number of incoming (resp., outgoing) edges attached to the skeleton vertex $r$ {\em in the skeleton graph}\, $\Ga_{sk}$ ({\em not}\, in the original graph $\Ga^\fw$).

\sip

Consider a filtration of the complex $\WGC_d$ by the number of skeleton vertices (which can not decrease under the action of the differential), and let $\{E_k\WGC_d, d_k\}_{k\geq 0}$ be the associated spectral sequence. By Maschke Theorem, we can assume without loss of generality that vertices and edges of the skeleton graph in the next step are distinguished so that the automorphism group of the sketelon graph is trivial.

\mip

{\sf Step 2: Computation of $E_1\WGC_d=H^\bu(E_0\WGC_d,d_0)$}. Let $\mathsf{tGC}_d\subset \dfcGC_d$ be the subspace
spanned by directed trees  of the form
\Beq\label{4: trees}
T_r=
\Ba{c}\resizebox{25mm}{!}{
\xy
 (0,10)*{^r},
 (0,8)*{\bu}="a",
(-10,0)*{\bu}="b_1",
(-2,0)*{\bu}="b_2",
(12,0)*{\bu}="b_3",
(2,-8)*{\bu}="c_1",
(-7,-8)*{\bu}="c_2",
(8,-8)*{\bu}="c_3",
(14,-8)*{\bu}="c_4",
(20,-8)*{\bu}="c_5",
(-11,-15)*{\bu}="d_1",
(-3,-15)*{\bu}="d_2",
\ar @{<-} "a";"b_1" <0pt>
\ar @{<-} "a";"b_2" <0pt>
\ar @{<-} "a";"b_3" <0pt>
\ar @{<-} "b_2";"c_1" <0pt>
\ar @{<-} "b_2";"c_2" <0pt>
\ar @{<-} "b_3";"c_3" <0pt>
\ar @{<-} "b_3";"c_4" <0pt>
\ar @{<-} "b_3";"c_5" <0pt>
\ar @{<-} "c_2";"d_1" <0pt>
\ar @{<-} "c_2";"d_2" <0pt>
\endxy
}\Ea
\Eeq
that is, by genus zero directed graphs whose every vertex (except the special one called {\em the root vertex}\, and denoted by $r$) has precisely one outgoing edge; the root vertex has no outgoing edges. Note that $T_r$ may consist solely of the root vertex, i.e.\ in general $\#V(T_r)\geq 1$.

\sip

The complex $E_0\WGC_d$ decomposes, up to the trivial ``skeleton" complex, into a product
$$
\prod_{\Ga_{sk}\in \dfcGC_d'} \bigotimes_{r\in V(\Ga_{sk})} C_r
$$
where
$$
C_r:= \prod_{T_r\in \mathsf{tGC}_d}\prod_{\fw: V(T) \rar \N_{\geq 1}} T_r^\fw
$$
where the second product is taken over all weight functions $\fw$ satisfying the standard condition
$\fw(v)\geq \max\{1, |v|_{in}-1\}$ for every non-root vertex $v\neq r$ of a tree $T$, and
$$
\fw(r)\geq \max\{1, |r|^{T}_{in} +  |r|^{sk}_{in} -1,  |r|^{sk}_{out}-1\}
$$
where  $ |r|^{T}_{in}$ (resp., $|r|^{sk}_{in}$) counts the number of incoming edges to $r$ in the tree $T$ (resp., in the skeleton graph $\Ga_{sk}$).

\sip

The induced differential is given by
\Beq\label{5: d_0 in tGC}
 d_0 T^\fw =\delta_0 T^\fw
\pm  \sum_{ v\in V(\Ga)\atop
 \fw(v)> \# In_v-1}\sum_{c\geq 1}
\Ba{c}
{ \xy
(0,0)*+{_c}*\cir{}="u1",
(0,8)*+{T^\fw}="ul",
\ar @{<-} "ul";"u1" <0pt>
\endxy}
\Ea
 \Eeq
where
\Bi
\item the first term $\delta_0 \Ga^\fw$ is equal to the sum $\pm\sum_{v\in V(\Ga)}\delta_{0,v}T^\fw$ over all vertices of $T$ (including the root vertex $r$) with $\delta_{0,v}\Ga^\fw$ obtained
from $T^\fw$ by replacing
   the vertex $v$  with weight, say, $a\in \N_{\geq 1}$ by the sum of graphs
   $$
 \xy
(0,0)*+{_{a}}*\cir{}="u1",
\endxy \ \ \lon \ \ \ \   \sum_{a=a'+ a''}
\Ba{c}   \resizebox{4mm}{!}
   {\xy
(0,10)*+{_{a'}}*\cir{}="u1",
(0,0)*+{_{a''}}*\cir{}="ul",
\ar @{->} "ul";"u1" <0pt>
\endxy}\Ea
   $$
   and taking one more sum over all possible way to reattach  the dangling  edges to the newly created vertices $\Ba{c}   \resizebox{4mm}{!}
   {\xy
(0,1)*+{_{a'}}*\cir{},
\endxy}\Ea$ and $\Ba{c}   \resizebox{4mm}{!}
   {\xy
(0,1)*+{_{a''}}*\cir{},
\endxy}\Ea$ in such a way
   that the condition (\ref{3: condition on weights}) is respected and no new skeleton vertices appear.

\item the second summand is a linear combination of trees obtained by attaching an in-vertex $\Ba{c}\resizebox{7mm}{!}{\xy
(0,1)*+{_c}*\cir{}="0",
(6,1)*{}="r",
\ar @{->} "0";"r" <0pt>
\endxy}\Ea$
 to each vertex $v$ of $T^\fw$ whose weight satisfies the strict inequality
$\fw(v)> |v|^{T}_{in}-1$ if $v\neq r$ and $\fw(r)> |v|^{T}_{in} + |v|^{sk}_{in} -1$ if $v$ is the root vertex $r$.

\Ei

Each complex $C_r$ has at least one cohomology class spanned by the root vertex assigned the following sum of weights,
 \Beq\label{5: weiht sum cohomology class}
O_r:=\sum_{a\geq \fs_r^{sk}} \ \xy(0,0)*+{_{a}}*\cir{};
\endxy
\Eeq
If the root vertex $r$ has skeleton valency satisfying  $|r|_{out}^{sk}\leq  |r|_{in}^{sk}$ and $|r|_{in}^{sk}\geq 2$,
then there is one more cohomology class
$$
o_r:= \xy(0,0)*+{_{\fs_r^{sk}}}*\cir{};\endxy
$$
spanned by the root vertex equipped with the minimal possible weight, that is with the skeleton weight.
Let us show that there are no other cohomology classes in the complexes $C_r$. Call a maximal connected subgraph of a tree $T_r^\fw$ a {\em string}\, if it consists of passing vertices and a univalent in-vertex (so the root vertex is never a part of some string); every string looks like an element (\ref{5: element of C}) of the auxiliary complex $C$ for some $n\geq 1$. If $\#V(T^\fw)\geq 2$, then $T^\fw$ contains at least one string. Call vertices of $T^\fw$ which belong to a string {\em stringy vertices}\, and the remaining vertices {\em core}\, ones. The set of core vertices is non-empty as it always contains the root vertex. Consider a filtration of $C_r$ by the number of core vertices, and then a filtration by the total weight of the stringy vertices. The induced differential acts only on stringy vertices in exactly the same way as on the vertices of the auxiliary complex $C$ (see (\ref{d on C}) above)  so that, by Lemma {\ref{2: lemma on C}}, each string (if any) has length 1 and is spanned by the single in-vertex of weight 1. If the number of core vertices is greater than or equal to $2$, then at least one core vertex contains at least two strings and hence does not contribute to the cohomology as $H^\bu(\odot^{\geq 2}C)=0$. Hence it remains to consider the case when the number of core vertices is equal to 1, i.e. the only core vertex is the root vertex. The next page of the spectral sequence by the total weight filtration of stringy vertices gives us a complex spanned by elements of two types --- by the single root vertex (equipped with an arbitrary weight $a\geq \fs_r^{sk}$) and by the root vertex together with one univalent in-vertex of weight 1 attached,
$$
\text{span} \left\langle
   \Ba{c}
\resizebox{5mm}{!}
{ \xy
(0,1)*+{_{a}}*\cir{}="2",
\endxy}
\Ea \ ,
 \Ba{c}
\resizebox{5mm}{!}
{ \xy
(0,-7)*+{_{1}}*\cir{}="1",
(0,1)*+{_{b}}*\cir{}="2",
\ar @{->} "1";"2" <0pt>
\endxy}
\Ea
\right\rangle
$$
where
$$
b\geq \left\{\Ba{ll}
\fs_r^{sk} & \text{if}\ |r|^{sk}_{out} > |r|^{sk}_{in} \ \text{or}\ \left( |r|^{sk}_{in} =1,  |r|^{sk}_{out}\leq 1\right) \\
\fs_r^{sk}+1 & \text{if}\ |r|^{sk}_{out} < |r|^{sk}_{in} \ \text{or}\ |r|^{sk}_{in} = |r|^{sk}_{out}\geq 2
\\
\Ea
\right.
$$
The induced differential acts trivially on the second generator and it acts on the first generator as follows
$$
 \Ba{c}
\resizebox{5mm}{!}
{ \xy
(0,1)*+{_{a}}*\cir{}="2",
\endxy}
\Ea \lon \left\{\Ba{cl}
 -\Ba{c}
\resizebox{12mm}{!}
{ \xy
(9,1)*+{_{1}}*\cir{}="1",
(0,1)*+{_{a}}*\cir{}="2",
\ar @{->} "1";"2" <0pt>
\endxy}
\Ea
& \text{if}\ \left(|r|^{sk}_{out} > |r|^{sk}_{in} \ \text{or}\ \left( |r|^{sk}_{in} =1,  |r|^{sk}_{out}\leq 1\right)\right)\ \text{and}\ a=\fs_r^{sk}\\
  \Ba{c}
\resizebox{12mm}{!}
{ \xy
(9,1)*+{_{1}}*\cir{}="1",
(0,1)*+{_{a-1}}*\cir{}="2",
\ar @{->} "1";"2" <0pt>
\endxy}
\Ea
-
\Ba{c}
\resizebox{12mm}{!}
{ \xy
(11,1)*+{_{1}}*\cir{}="1",
(0,1)*+{_{a}}*\cir{}="2",
\ar @{->} "1";"2" <0pt>
\endxy}
\Ea
& \text{if}\ |r|^{sk}_{out} > |r|^{sk}_{in} \ \text{or}\ |r|^{sk}_{in} = |r|^{sk}_{out}=1\ a>\fs_r^{sk}\\
0 & \text{if}\  |r|^{sk}_{out} < |r|^{sk}_{in} \ \text{or}\ |r|^{sk}_{in} = |r|^{sk}_{out}\geq 2\ \text{and}\ a=\fs_r^{sk}\\
 -\Ba{c}
\resizebox{12mm}{!}
{ \xy
(9,1)*+{_{1}}*\cir{}="1",
(0,1)*+{_{a}}*\cir{}="2",
\ar @{->} "1";"2" <0pt>
\endxy}
\Ea& \text{if}\  |r|^{sk}_{out} < |r|^{sk}_{in} \ \text{or}\ |r|^{sk}_{in} = |r|^{sk}_{out}\geq 2\ \text{and}\ a=\fs_r^{sk}+1\\
  \Ba{c}
\resizebox{12mm}{!}
{ \xy
(9,1)*+{_{1}}*\cir{}="1",
(0,1)*+{_{a-1}}*\cir{}="2",
\ar @{->} "1";"2" <0pt>
\endxy}
\Ea
-
\Ba{c}
\resizebox{12mm}{!}
{ \xy
(11,1)*+{_{1}}*\cir{}="1",
(0,1)*+{_{a}}*\cir{}="2",
\ar @{->} "1";"2" <0pt>
\endxy}
\Ea
& \text{if}\  |r|^{sk}_{out} < |r|^{sk}_{in} \ \text{or}\ |r|^{sk}_{in} = |r|^{sk}_{out}\geq 2\ \text{and}\ a>\fs_r^{sk}+1\\
\Ea
\right.
$$
 It is now almost immediate to see that the cohomology is generated by the class $O_r$ and, in the case when $|r|_{out}^{sk}\leq  |r|_{in}^{sk}$ and $|r|_{in}^{sk}\geq 2$,
 by one extra the class $o_r$, i.e. that it is at most 2-dimensional.

\sip

We conclude that the next page $E_1\WGC_d=H^\bu(E_0\WGC_d,d_0)$ is generated by weighted graphs $\Ga^\fw$ with no univalent in-vertices and with every vertex $v$ either assigned the sum of weights as in (\ref{5: weiht sum cohomology class}) (we call such vertices the {\em hairy}\, ones and denote them by double circles $\Ba{c}\resizebox{5mm}{!}
{ \xy
(0,1)*+{_{\ }}*\cir{},
(0,1)*+{ \ \, }*\cir{},
\endxy}
\Ea$) or, if  $|v|_{out}\leq  |v|_{in}$ and $|v|_{in}\geq 2$,
the weight $\fs_v^{sk}$; we call the latter vertices the {\em bold}\, ones.

\mip

{\sf Step 3: Study of the complex $(E_1\WGC_d, d_1)$}.  The induced differential $d_1$ is given by
\Beq\label{5: d_1 in E_1WGC}
 d_1\Ga =\sum_{v\in V(\Ga)} d_1^v \Ga
\pm  \sum_{ v\ \text{is baldand}\atop
 |r|^{sk}_{out} < |r|^{sk}_{in}}
\Ba{c}
{ \xy
(0,0)*+{_{\ }}*\cir{},
(0,0)*+{ \ \ }*\cir{}="u1",
(0,-8)*+{\Ga}="ul",
\ar @{->} "ul";"u1" <0pt>
\endxy}
\Ea
\pm  \sum_{ v\ \text{is hairy}}
\Ba{c}
{ \xy
(0,0)*+{_{\ }}*\cir{},
(0,0)*+{ \ \ }*\cir{}="u1",
(0,-8)*+{\Ga}="ul",
\ar @{->} "ul";"u1" <0pt>
\endxy}
\Ea
 \Eeq
 where the operator $d_1^v$ splits a bald vertex $v$ into a pair bald vertices as in (\ref{5: splitting weight vertex}) and a hairy vertex $v$ into a pair of hairy vertices as follows
$$
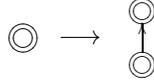
\Ba{c}\resizebox{5mm}{!}
{ \xy
(0,1)*+{_{\ }}*\cir{},
(0,1)*+{ \ \ }*\cir{},
\endxy}
\Ea
\lon
\Ba{c}
 \resizebox{4.5mm}{!}
{ \xy
(0,-4.5)*+{_{\ }}*\cir{},
(0,-4.5)*+{ \ \ }*\cir{}="0",
(0,4.5)*+{_{\ }}*\cir{},
(0,4.5)*+{\ \ }*\cir{}="1",
\ar @{->} "0";"1" <0pt>
\endxy}
\Ea
$$
Two other summands stand for the sums over attaching a new hairy univalent-out vertex to bald vertices
$v$ with  $|r|^{sk}_{out} < |r|^{sk}_{in}$ and to hairy vertices.

\sip

Note that passing and univalent (out-)vertices (if any) in a graph $\Ga$ from $E_1\WGC_d$ must be hairy.

\sip

Deleting from a graph $\Ga\in E_1\WGC_d^\star$ recursively all univalent (hairy) out-vertices, we obtain a graph $\Ga'$ called the {\em prime graph}\, of $\Ga$ (cf.\ Step 1 above). Consider next a filtration by the number of vertices in the prime graphs, and let $(E_{r,1}\WGC_d, d_{r,1})$ be the associated spectral sequence. The associated graded complex $(E_{0,1}\WGC_d, d_{0,1})$ is the tensor product of tree type complexes spanned by elements of the form (\ref{4: trees}) with arrows reversed on every edge, and with every vertex except the root being a hairy vertex. The root vertex $r$ can be a hairy vertex or a bald vertex satisfying the conditions  $|r|^{prime}_{out} < |r|^{prime}_{in}$ and
$|r|^{prime}_{in}\geq 2$; here $|r|^{prime}_{in}$ and $|r|^{prime}_{out}$ stand for the valencies
of the root vertex in the {\em prime graph}\, $\Ga'$ (not in the original graph $\Ga$).
An analysis similar to the one discussed above leads us easily to the conclusion that the next page
$$
E_{1,1}\WGC_d\simeq H^\bu((E_{0,1}\WGC_d, d_{0,1})
$$
 is spanned by graphs with no univalent vertices and with no bald vertices $v$ such that $|v|_{in}\neq |v|_{out}$. The induced differential $d_{1,1}$ splits a bald vertex into the pair of bald vertices,  and a hairy vertex into a pair of hairy vertices.
This complex decomposes naturally into a direct sum,
$$
E_{1,1}\WGC_d= E_{1,1}\WGC_d^{bold} \oplus E_{1,1}\WGC_d^{hairy} \oplus E_{1,1}\WGC_d^{mixed}
$$
where
\Bi
\item the summand $E_{1,1}\WGC_d^{bold}$ is spanned by graphs whose every vertex $v$ is bald and hence satisfies $|v|_{in}=|v|_{out}$. This subcomplex is precisely the summand $\WGC_d^{=}\simeq \Der^=(\wHoQPcd)$ in the decomposition (\ref{5: Decomp od WGC^star}).
\item The summand $E_{1,1}\WGC_d^{hairy}$ is exactly the image, $\cF^\circlearrowright(\dfcGC_d^{\geq 2})$, of the injective morphism  of complexes (\ref{4: F from dcGC to wGC}).

\item The summand $ E_{1,1}\WGC_d^{mixed}$ is spanned by graphs with at least one vertex bald and at least one hairy vertex.

\Ei
Hence to complete the proof of the Theorem it remains to show that the subcomplex $ E_{1,1}\WGC_d^{mixed}$ is acyclic. For this let us call a maximal connected path of consisting of passing vertices (if any) a {\em weighted edge} (recall that every passing vertex must be a hairy vertex).  One can understand such a weighted edge as an ordinary edge together with a positive number (its {\em weight}) assigned to it,
$$
\stackrel{k}{\lon}
$$
where $k$ counts the number of passing vertices, e.g.
$$
\stackrel{2}{\lon}\, :=\Ba{c}
 \resizebox{18mm}{!}
{ \xy
(-39,1)*{}="-1",
(-30,1)*+{_{\ }}*\cir{},
(-30,1)*+{ \ \ }*\cir{}="0",
(-20,1)*+{_{\ }}*\cir{},
(-20,1)*+{\ \ }*\cir{}="1",
(-11,1)*{}="2",
\ar @{->} "-1";"0" <0pt>
\ar @{->} "0";"1" <0pt>
\ar @{->} "1";"2" <0pt>
\endxy}
\Ea
$$
We shall allow next weighted edges equipped with the zero weight by identifying them with ordinary edges,
$$
\stackrel{0}{\lon} \, :=\,  \lon.
$$

Any graph $\Ga$ from $ E_1\WGC_d^{mixed}$ can be understood as a graph $\Ga^{reduced}$ with no passing vertices but with weighted edges. Every such a graph containing at least one weighted wedge connecting bald and hairy vertices; the number of such edges is preserved under the differential
$d_{1,1}$ so that we can assume without loss of generality that such {\em mixed}\, edges are distinguished, e.g.\ marked by positive integers starting with 1.
Consider a filtration by the number of vertices in $\Ga^{reduced}$ plus the total weight of mixed edges with labels $\geq 2$. The induced differential acts only on the mixed edge labelled by $1$ increasing its weight from $k\geq 0$ to $k+1$.  Hence the complex  $ E_1\WGC_d^{mixed}$ is acyclic indeed, and the Theorem is proven.
\end{proof}




The above Theorem, the isomorphism
$$
\Der(\wHoQP_{c,d})=\odot^{\bu\geq 1}(\WGC_{c+d+1}[-1-c-d])[1+c+d]
$$
and the canonical quasi-isomorphism (\ref{2: from GC to dcGC}) imply immediately Theorem {\ref{1: Theorem on FGC and Der}} in the Introduction.

\sip

An immediate corollary to Theorem {\ref{1: Theorem on FGC and Der}} is that the canonical morphism of dg Lie algebras
$$
\Der(\HoLB_{c,d}^{\atop \bigstar \circlearrowright}) \lon \Der(\HoQP_{c,d}^{\atop \bigstar \circlearrowright})
$$
is a quasi-isomorphism.

\sip

Considering an oriented version of the weighted graph complex $\WGC_d$ and using arguments almost identical to the ones used in the proof of Theorem {\ref{4: Theorem on dcGC and dwGC}} leads us to the the following result about the derivation complex of the ordinary properad of $\Z$-graded Poisson structures.

\subsubsection{\bf Theorem}
{\em There is a canonical morphism of dg Lie algebras,
$$
\mathsf{OGC}^{\geq 2}_{c+d+1} \lon \mathsf{Der}(\HoQP_{c,d})
$$
which is a quasi-isomorphism up to one rescaling class $r$. In particular, there is an isomorphism of Lie algebras
$$
H^0(\Der(\HoQP_{0,1}))=\K
$$
that is, the only homotopy non-trivial automorphism of the properad $\widehat{\HoQP}_{0,1}$ is given by the rescaling of the generators as follows
$$
\Ba{c}\resizebox{14mm}{!}{\begin{xy}
 <0mm,0mm>*{\circ};<0mm,0mm>*{}**@{},
 <-0.6mm,0.44mm>*{};<-8mm,5mm>*{}**@{-},
 <-0.4mm,0.7mm>*{};<-4.5mm,5mm>*{}**@{-},
 <0mm,0mm>*{};<-1mm,5mm>*{\ldots}**@{},
 <0.4mm,0.7mm>*{};<4.5mm,5mm>*{}**@{-},
 <0.6mm,0.44mm>*{};<8mm,5mm>*{}**@{-},
   <0mm,0mm>*{};<-8.5mm,5.5mm>*{^1}**@{},
   <0mm,0mm>*{};<-5mm,5.5mm>*{^2}**@{},
   <0mm,0mm>*{};<4.5mm,5.5mm>*{^{n\hspace{-0.5mm}-\hspace{-0.5mm}1}}**@{},
   <0mm,0mm>*{};<9.0mm,5.5mm>*{^n}**@{},
 <-0.6mm,-0.44
 mm>*{};<-8mm,-5mm>*{}**@{-},
 <-0.4mm,-0.7mm>*{};<-4.5mm,-5mm>*{}**@{-},
 <0mm,0mm>*{};<-1mm,-5mm>*{\ldots}**@{},
 <0.4mm,-0.7mm>*{};<4.5mm,-5mm>*{}**@{-},
 <0.6mm,-0.44mm>*{};<8mm,-5mm>*{}**@{-},
   <0mm,0mm>*{};<-8.5mm,-6.9mm>*{^1}**@{},
   <0mm,0mm>*{};<-5mm,-6.9mm>*{^2}**@{},
   <0mm,0mm>*{};<4.5mm,-6.9mm>*{^{n\hspace{-0.5mm}-\hspace{-0.5mm}1}}**@{},
   <0mm,0mm>*{};<9.0mm,-6.9mm>*{^n}**@{},
 \end{xy}}\Ea
 \lon
 \la^{n-1}
 \Ba{c}\resizebox{14mm}{!}{\begin{xy}
 <0mm,0mm>*{\circ};<0mm,0mm>*{}**@{},
 <-0.6mm,0.44mm>*{};<-8mm,5mm>*{}**@{-},
 <-0.4mm,0.7mm>*{};<-4.5mm,5mm>*{}**@{-},
 <0mm,0mm>*{};<-1mm,5mm>*{\ldots}**@{},
 <0.4mm,0.7mm>*{};<4.5mm,5mm>*{}**@{-},
 <0.6mm,0.44mm>*{};<8mm,5mm>*{}**@{-},
   <0mm,0mm>*{};<-8.5mm,5.5mm>*{^1}**@{},
   <0mm,0mm>*{};<-5mm,5.5mm>*{^2}**@{},
   <0mm,0mm>*{};<4.5mm,5.5mm>*{^{n\hspace{-0.5mm}-\hspace{-0.5mm}1}}**@{},
   <0mm,0mm>*{};<9.0mm,5.5mm>*{^n}**@{},
 <-0.6mm,-0.44
 mm>*{};<-8mm,-5mm>*{}**@{-},
 <-0.4mm,-0.7mm>*{};<-4.5mm,-5mm>*{}**@{-},
 <0mm,0mm>*{};<-1mm,-5mm>*{\ldots}**@{},
 <0.4mm,-0.7mm>*{};<4.5mm,-5mm>*{}**@{-},
 <0.6mm,-0.44mm>*{};<8mm,-5mm>*{}**@{-},
   <0mm,0mm>*{};<-8.5mm,-6.9mm>*{^1}**@{},
   <0mm,0mm>*{};<-5mm,-6.9mm>*{^2}**@{},
   <0mm,0mm>*{};<4.5mm,-6.9mm>*{^{n\hspace{-0.5mm}-\hspace{-0.5mm}1}}**@{},
   <0mm,0mm>*{};<9.0mm,-6.9mm>*{^n}**@{},
 \end{xy}}\Ea, \ \ \ \forall \la \in \K\setminus 0, \ \forall\ n\geq 2.
 $$
}

\sip

The second claim follows from the well-known equality of cohomology groups \cite{Wi2}
$$
H^\bu(\OGC_2)=H^\bu(\GC_1^{\geq 2})
$$
and the fact that $H^0(\GC_1^{\geq 2})=0$.

\bip


{\Large
\section{\bf Classification of universal quantizations of quadratic  Poisson structures}
}

\mip


\subsection{Kontsevich formality map as a morphism of dg operads}
Let $V$ be a finite-dimensional $\Z$-graded vector space, and $V^*:=\Hom(V,\K)$ be its dual vector space. Let us view $V^*$ as a  pointed manifold and view the symmetric tensor algebra
$$
\f_V=\bigoplus_{k\geq 0}  \f_V^k, \ \ \  \f_V^k:=\odot^k V, 
$$
as the algebra of formal smooth functions on $V^*$. Deformations of the standard graded commutative product in $\f_V$ are controlled by the standard dg Lie algebra, the so called {\em Hochschild complex}, $$
\mathsf{Hoch}(\f_V,\f_V):= \bigoplus_{n\geq 0} \mathsf{Hoch}^n(\f_V,\f_V), \ \ \ \   \mathsf{Hoch}^n(\f_V,\f_V):=\Hom(\ot^n \f_V, \f_V)[1-n]
$$
of polydifferential operators (see \cite{Ko2} for the explicit formulae of the differential in $\mathsf{Hoch}(\f_V,\f_V)$ and the Lie brackets).

\sip

The space of tangent vector fields on $V^*$
can be identified with the Lie algebra of derivations of $\f_V$ which are uniquely determined by its values on the generators, that is by the values on the elements of $V$,
$$
\Der(\f_V)\simeq \Hom(V,{\odot^\bu} V)
$$
The Lie algebra of polyvector fields on $V^*$ can be identified as a vector space with
$$
 \sT(V):= \wedge^\bu_{\f_V} \Der(\f_V)  \simeq \bigoplus_{k\geq 0} \sT^k(V), \ \ \
 \sT^k(V):=\bigoplus_{l}
 \Hom( \wedge^l V,\odot^k V)[1-l]
$$
The explicit formula for the Lie bracket in $\sT(V)$ can be found, for example, in \cite{Ko2}. The differential in $\sT(V)$ is assumed to be zero.

\sip

Let $\hbar$ be a formal parameter (of cohomological degree zero), and let, for a vector space $A$,
the symbol $A[[\hbar]]$ stand for the vector space of formal power series in $\hbar$ with coefficients in $A$. Maxim Kontsevich formality map is a quasi-isomorphism of dg Lie algebras
\Beq\label{5: full formality map}
\sF^K: \sT(V)[[\hbar]]\lon \mathsf{Hoch}(\f_V,\f_V)[[\hbar]],
\Eeq
An explicit formula for such a quasi-isomorphism for $\K=\R$ was constructed in \cite{Ko2} (implying its existence over any field $\K$ of characteristic zero). That formula gives a {\em universal}\, formality morphism which is applicable to an arbitrary finite-dimensional  graded vector space $V$, and hence can be reformulated without any reference to any particular $V$. Such a reformulation was given in \cite{AM} as a morphism of dg operads
\Beq\label{5: quant map F}
F: c\cA ss_\infty \lon \f(\swHoLB_{0,1})
\Eeq
satisfying a certain non-triviality condition.
Here $c\cA ss_\infty$ is a  dg operad of {\em curved
$A_\infty$-algebras} which is, by definition,  the free operad generated by the $\bS$-module
\Beq\label{5: generators of cAssinfty}
E(n):=\K[\bS_n][n-2]= \mbox{span}
\left(\Ba{c}\resizebox{20mm}{!}{
\begin{xy}
 <0mm,0mm>*{\bullet};<0mm,0mm>*{}**@{},
 <0mm,0mm>*{};<-8mm,-5mm>*{}**@{-},
 <0mm,0mm>*{};<-4.5mm,-5mm>*{}**@{-},
 <0mm,0mm>*{};<0mm,-4mm>*{\ldots}**@{},
 <0mm,0mm>*{};<4.5mm,-5mm>*{}**@{-},
 <0mm,0mm>*{};<8mm,-5mm>*{}**@{-},
   <0mm,0mm>*{};<-11mm,-7.9mm>*{^{\sigma(1)}}**@{},
   <0mm,0mm>*{};<-4mm,-7.9mm>*{^{\sigma(2)}}**@{},
   <0mm,0mm>*{};<10.0mm,-7.9mm>*{^{\sigma(n)}}**@{},
 <0mm,0mm>*{};<0mm,5mm>*{}**@{-},
 \end{xy}}\Ea\ \
\right)_{\sigma\in \bS_n},\ \ \ \forall\ n\geq 0
\Eeq
where $\K[\bS_n]$ is the group algebra of the permutation group. It has a canonical basis is given by permutations $\sigma\in \bS_n$ which we identify with planar corollas with one output leg and $n$ input legs (every such a corolla stands, roughly speaking, for a generic lineal map $\ot^n V \lon V[2-n]$).
The differential in $c\cA ss_\infty$ is given on the generators by the formula
$$
\delta
\Ba{c}\resizebox{10mm}{!}{
\begin{xy}
<0mm,0mm>*{\bullet},
<0mm,5mm>*{}**@{-},
<-5mm,-5mm>*{}**@{-},
<-2mm,-5mm>*{}**@{-},
<2mm,-5mm>*{}**@{-},
<5mm,-5mm>*{}**@{-},
<0mm,-7mm>*{_{1\ \ \ \ldots\ \ \ n}},
\end{xy}}\Ea
=\sum_{k=0}^{n}\sum_{l=0}^{n-k}
(-1)^{k+l(n-k-l)+1}
\Ba{c}\resizebox{30mm}{!}{
\begin{xy}
<0mm,0mm>*{\bullet},
<0mm,5mm>*{}**@{-},
<4mm,-7mm>*{^{1\ \ \dots \ \ k\qquad\ \    (k+l+1)\ \ \dots\ n}},
<-14mm,-5mm>*{}**@{-},
<-6mm,-5mm>*{}**@{-},
<20mm,-5mm>*{}**@{-},
<8mm,-5mm>*{}**@{-},
<0mm,-5mm>*{}**@{-},
<0mm,-5mm>*{\bullet};
<-5mm,-10mm>*{}**@{-},
<-2mm,-10mm>*{}**@{-},
<2mm,-10mm>*{}**@{-},
<5mm,-10mm>*{}**@{-},
<0mm,-12mm>*{_{k+1\ \dots\ k+l }},
\end{xy}}\Ea.
$$

\sip

The properad $\swHoLB_{0,1}$ has been defined in \S 2. The symbol $\f$ stands for a polydifferential functor
$$
\f: \text{\sf Category of dg (wheeled) props} \lon \text{\sf Category of dg operads}
$$
introduced in \cite{MW2}. Elements of $\f(\swHoLB_{0,1})$ are generated by graphs from $\swHoLB_{0,1}$
whose outgoing legs are symmetrized and attached to the new white out-vertex, while ingoing legs are partitions into disjoint union of some subsets, legs in each subset are symmetrized and are attached to a new in-vertex which is labelled by an integer. For example \cite{AM}, For example, an element
$
e=\vspace{-1mm}\Ba{c}\resizebox{11mm}{!}{
\xy
(0,0)*{\circ}="o",
(-2,5)*{}="2",
(4,5)*{\circ}="3",
(4,10)*{}="u",
(4,0)*{}="d1",
(7,0)*{}="d2",
(10,0)*{}="d3",
(-1.5,-5)*{}="5",
(1.5,-5)*{}="6",
(4,-5)*{}="7",
(-4,-5)*{}="8",
(-2,7)*{_1},
(4,12)*{_2},
(-1.5,-7)*{_2},
(1.5,-7)*{_3},
(10.4,-1.6)*{_6},
(-4,-7)*{_1},
(4,-1.6)*{_4},
(7,-1.6)*{_5},
\ar @{-} "o";"2" <0pt>
\ar @{-} "o";"3" <0pt>
\ar @{-} "o";"5" <0pt>
\ar @{-} "o";"6" <0pt>
\ar @{-} "o";"8" <0pt>
\ar @{-} "3";"u" <0pt>
\ar @{-} "3";"d1" <0pt>
\ar @{-} "3";"d2" <0pt>
\ar @{-} "3";"d3" <0pt>
\endxy}\Ea
 \in \widehat{\HoLB}_{0,1}^{\atop \bigstar\circlearrowright}(2,6)
$
can generate the following element
$
\Ba{c}\resizebox{17mm}{!}{ \xy
(-1.5,5)*{}="1",
(1.5,5)*{}="2",
(9,5)*{}="3",
 (0,0)*{\circ}="A";
  (9,3)*{\circ}="O";
   (5,12)*+{\hspace{2mm}}*\frm{o}="X";
 (-6,-10)*+{_1}*\frm{o}="B";
  (6,-10)*+{_2}*\frm{o}="C";
   (14,-10)*+{_3}*\frm{o}="D";
    (22,-10)*+{_4}*\frm{o}="E";
 "A"; "B" **\crv{(-5,-0)}; 
  "A"; "D" **\crv{(5,-0.5)};
  "A"; "C" **\crv{(-5,-7)};
   "A"; "O" **\crv{(5,5)};
\ar @{-} "O";"C" <0pt>
\ar @{-} "O";"D" <0pt>
\ar @{-} "O";"X" <0pt>
\ar @{-} "A";"X" <0pt>
\ar @{-} "O";"B" <0pt>
 \endxy}
 \Ea$ in the operad $\f(\widehat{\HoLB}_{0,1}^{\atop \bigstar\circlearrowright})(4)$.

 \mip


It was shown in \cite{AM} that for any universal formality morphism (\ref{5: quant map F})
there is a canonically associated morphism
 of complexes
 $$
  \mathsf{FGC}^{\geq 2}_2 \lon \Def\left(c\Ass_\infty \stackrel{F}{\lon} \f(\wh{\HoLB}_{0,1}^{\atop \bigstar\circlearrowright})\right)[1]
 $$
 which is a quasi-isomorphism, where the symbol $\Def$ stands for the standard deformation complex of the particular morphism $\cF$ constructed along the recipe given in \cite{MV} and its slight modification introduced in \S 3.4 of \cite{MW3}. This results gives us (almost immediately) the classification
 of all universal formality morphisms up to homotopy equivalence --- the set of such morphisms can be identified with set of Drinfeld associators. In the next subsection we obtain a similar conclusion for the set of homotopy classes of {\em homogeneous}\, formality morphisms.

\subsection{Kontsevich formality map applied to $\Z$-graded quadratic Poisson structures}
 Let
$$
\mathsf{Hoch}_{(0)}(\f_V,\f_V)= \bigoplus_{n\geq 0}\mathsf{Hoch}^n_{(0)}(\f_V,\f_V), \ \ \
\mathsf{Hoch}^n_{(0)}(\f_V,\f_V):=\Hom(\f_V^{k_1}\ot\ldots \ot \f_V^{k_n}, \f_V^{k_1+...+k_n})[1-n]
$$
be the subspace of $\mathsf{Hoch}(\f_V,\f_V)$ spanned by homogeneous polydifferential operators, that is, the ones which preserve the total polynomial degree of formal functions. It is a dg Lie subalgebra.

\sip

 Similarly, let
$$
 \sT_{(0)}(V):= \wedge \Der_{(0)}(\f_V)  \simeq \bigoplus_{k}\Hom(\wedge^k V,\odot^k V)[1-k]
$$
be the Lie subalgebra of $\sT(V)$ spanned by homogeneous polyvector fields.
  The Maurer-Cartan elements
of the Lie algebra $\sT_{(0)}(V)$ are called {\em $\Z$-graded quadratic Poisson structures}. They have a decomposition,
\Beq\label{1: pi formal power series}
\pi=\sum_{n=0}^\infty\pi_n,\ \ \ \  \ \pi_n\in Hom(\wedge^n V, \odot^n V)[2-n]
\Eeq
If the vector space $V$ is concentrated in degree zero, then only the term with $n=2$ survives giving the notion of ordinary quadratic Poisson structure. The term $\pi_0$ is just a constant playing no role in deformation quantization theory; hence it can be ignored.

\sip

Using properties P1-P5 (see \S 7 in \cite{Ko2}) of the Kontsevich formality $\cF^K$, it is easy to see that the restriction of the map $\sF^K: \sT(V)[[\hbar]]\lon \mathsf{Hoch}(\f_V,\f_V)[[\hbar]]$
to the subcomplex
$$
\sT(V)_{(0)}^{\geq 1}:=\bigoplus_{k\geq 1}\sT(V)_{(0)}^{k} \subset \sT(V)
$$
takes values in
$$
\mathsf{Hoch}_{(0)}^{\geq 1}(\f_V,\f_V)[[\hbar]]=\bigoplus_{k\geq 1}\mathsf{Hoch}_{(0)}^{k}(\f_V,\f_V)[[\hbar]].
$$
Put another way, the Kontsevich map gives us a quasi-isomorphism of reduced complexes,
\Beq\label{5: homog formality}
\sF^K: \sT(V)_{(0)}^{\geq 1} \lon \mathsf{Hoch}_{(0)}^{\geq 1}(\f_V,\f_V)[[\hbar]]
\Eeq
Note that such a reduction is not possible in the case (\ref{5: full formality map}), with
the summand $\mathsf{Hoch}^{0}(\f_V,\f_V)[[\hbar]]$ in the r.h.s.\ playing a key role of the classification of homotopy inequivalent formality maps \cite{Do,AM}. Let us call any strongly homotopy quasi-isomorphism of Lie algebras as in (\ref{5: homog formality}) a {\em homogeneous}\, formality map.

\sip

Our purpose is to classify all {\em universal}\, homotopy inequivalent homogeneous formality maps. To give a precise definition of what {\em universal}\, means we shall use again the theory of operads.
Let $\cA ss_\infty$ be the standard operad of (flat) strongly homotopy associative algebras. It is a free operad generators by the $\bS$ module (\ref{5: generators of cAssinfty}) with $n\geq 2$ only.
The canonical epimorphism of dg props
$$
\wh{\HoLB}_{0,1}^{\atop \bigstar\circlearrowright} \lon \wh{\HoQP}_{0,1}^{\atop\circlearrowright}
$$
induces an epimorphism of dg operads
$$
\f\left(\wh{\HoLB}_{0,1}^{\atop \bigstar\circlearrowright}\right) \lon \f\left(\wh{\HoQP}_{0,1}^{\atop \circlearrowright}\right).
$$
Let $\f^{\geq 1}\left(\wh{\HoQP}_{0,1}^{\atop \circlearrowright}\right)$ be a subspace of
$\f\left(\wh{\HoQP}_{0,1}^{\atop \circlearrowright}\right)$ spanned by graphs with at least one incoming white vertex.
The Kontsevich formality map (\ref{5: homog formality}) gives us a morphism of dg operads
\Beq\label{F from Assinfty to O(QP)}
\cF^K: \cA ss_\infty \lon \f^{\geq 1}\left(\wh{\HoQP}_{0,1}^{\atop \circlearrowright}\right)
\Eeq
satisfying the following boundary conditions, $O(2)$)
\Beq\label{5: Boundary cond for formality map}
\cF^K\left(
\Ba{c}\resizebox{10mm}{!}{
\begin{xy}
<0mm,0mm>*{\bullet},
<0mm,5mm>*{}**@{-},
<-5mm,-5mm>*{}**@{-},
<-2mm,-5mm>*{}**@{-},
<2mm,-5mm>*{}**@{-},
<5mm,-5mm>*{}**@{-},
<0mm,-7mm>*{_{1\ \ \ \ldots\ \ \ n}},
\end{xy}}\Ea
\right)=\left\{
\Ba{ll}
 \Ba{c}\resizebox{9mm}{!}{ \xy
   (0,9)*+{\hspace{2mm}}*\frm{o}="X";
 (-5,0)*+{_1}*\frm{o}="B";
  (5,0)*+{_2}*\frm{o}="C";
 \endxy}\Ea
 + \sum_{p\geq 0}\frac{1}{2}
 \Ba{c}\resizebox{10mm}{!}{ \xy
 (0,-3)*{\circ}="a";
   (0,9)*+{\hspace{2mm}}*\frm{o}="X";
 (-7,-9)*+{_1}*\frm{o}="B";
  (7,-9)*+{_2}*\frm{o}="C";
   "a"; "X" **\crv{(-5,-0)}; 
   "a"; "X" **\crv{(+5,-0)}; 
  \ar @{-} "a";"B" <0pt>
  \ar @{-} "a";"C" <0pt>
 \endxy}\Ea + O(2) & \text{if}\ n=2\\
\frac{{1}}{n!}
 \Ba{c}\resizebox{13mm}{!}{ \xy
 (0,2.8)*{^n};
 (0,1)*{...};
 (3.5,-10)*{...};
 (0,-3)*{\circ}="a";
   (0,9)*+{\hspace{2mm}}*\frm{o}="X";
 (-10,-10)*+{_1}*\frm{o}="B";
  (-3,-10)*+{_2}*\frm{o}="C";
  (10,-10)*+{_{n}}*\frm{o}="E";
   "a"; "X" **\crv{(-5,-0)}; 
   "a"; "X" **\crv{(+5,-0)}; 
   "a"; "X" **\crv{(9,-0)}; 
   "a"; "X" **\crv{(-9,-0)}; 
  \ar @{-} "a";"B" <0pt>
  \ar @{-} "a";"C" <0pt>
  \ar @{-} "a";"E" <0pt>
 \endxy}\Ea + O(2) & \text{for $n\geq 3$}\\
 \Ea
\right.
\Eeq
 where $O(2)$ stand for the terms with the number of internal vertices $\geq 2$. Any morphism $\cF_0$ of dg operads as in (\ref{F from Assinfty to O(QP)}) satisfying the above boundary condition (which guarantees its non-triviality) is called a {\em universal homogeneous formality map}.

\subsection{Theorem (Classification of homogeneous formality maps)}\label{5: Corollary on GCor and Def(assb to Dlie)} {\em For any homogeneous formality morphism
$$
\cF_0:
\Ass_\infty {\lon} \f^{\geq 1}(\wh{\HoQP}_{0,1}^{\atop \circlearrowright})
$$
there is a canonically associated morphism
 of complexes
 $$
f_0:  \mathsf{FGC}^{\geq 2}_2 \lon \Def\left(\Ass_\infty \stackrel{\cF_0}{\lon} \f^{\geq 1}(\wh{\HoQP}_{0,1}^{\atop \circlearrowright})\right)[1]
 $$
 which is a quasi-isomorphism.}
\begin{proof} The proof of this Theorem is a straightforward adoption of the arguments used in the proof of Proposition 5.4.1 in \cite{MW3}, and is based essentially on the contractibility of the permutahedra polytopes. Let us first explain the structure of the induced morphism $f_{0}$.
Any derivation of the dg wheeled prop $\wh{\HoQP}_{0,1}^{\atop \circlearrowright}$ gives us an infinitesimal deformation of the identity automorphism
of $\wh{\HoQP}_{0,1}^{\atop \circlearrowright}$ which in turn induces
an infinitesimal deformation of the identity automorphism of the dg operad $\f^{\geq 1}(\wh{\HoQP}_{0,1}^{\atop \circlearrowright})$ which in its turn
induces an infinitesimal deformation of the
homogeneous formality map $\cF$ via the action on the r.h.s. Therefore we have a sequence of canonical morphisms of complexes
\Beqrn
\Der(\wh{\HoQP}_{0,1}^{\atop \circlearrowright}) \stackrel{\simeq}{\lon} \Def\left(\wh{\HoQP}_{0,1}^{\atop \circlearrowright}\stackrel{\Id}{\rar} \wh{\HoQP}_{0,1}^{\atop \circlearrowright}\right)[1] \stackrel{\simeq}{\lon}   &\Def\left(\f^{\geq 1}(\wh{\HoLB}_{0,1}^{\atop \circlearrowright})\stackrel{\Id}{\rar}
\f^{\geq 1}(\wh{\HoLB}_{0,1}^{\atop \circlearrowright})\right)[1]& \\
& \downarrow &\\
& \Def\left(\Ass_\infty \stackrel{\cF_0}{\lon} \f^{\geq 1}(\wh{\HoQP}_{0,1}^{\atop \circlearrowright})\right)[1] &
\Eeqrn
The first arrow above is a quasi-isomorphism because both complexes are identical to each other up to a degree shift. The second arrow is a quasi-isomorphism because the functor $\f$ is exact \cite{MW2}. We shall show below that the third arrow is also a quasi-isomorphism. Combining this latest claim with the quasi-isomorphism (\ref{1: map from FGC to DerQP}) proves the Theorem.

\sip
 Thus to prove the Theorem it is enough to show that the composition
\Beq\label{5: g_cF}
c: \Der(\wh{\HoQP}_{0,1}^{\atop \circlearrowright}) \lon
\Def\left(\Ass_\infty \stackrel{\cF_0}{\lon} \f^{\geq 1}( \wh{\HoQP}_{0,1}^{\atop \circlearrowright})\right)[1]
\Eeq
is a quasi-isomorphism.
Both complexes in (\ref{5: g_cF}) admit filtrations by the number of edges which is preserved by the map $c$. Hence that map induces a morphism of the associated spectral sequences,
$$
c^r: (\cE_r\Der(\wh{\HoQP}_{0,1}^{\atop \circlearrowright}), d_r) \lon \left(\cE_r
\Def\left(\Ass_\infty \stackrel{\cF}{\lon} \f^{\geq 1}( \wh{\HoQP}_{0,1}^{\atop \circlearrowright})\right)[1], \delta_r\right).
$$
The induced differential $d_0$ on the initial page of the spectral sequence of the l.h.s.\ is trivial, $d_0=0$. The induced differential $\delta_0$ on the initial page of the spectral sequence of the r.h.s. is not trivial and is determined by the
following summand in $\cF_0$ (see the boundary condition (\ref{5: Boundary cond for formality map}) for $\cF_0$),
$$
 \Ba{c}\resizebox{9mm}{!}{ \xy
   (0,9)*+{\hspace{2mm}}*\frm{o}="X";
 (-5,0)*+{_1}*\frm{o}="B";
  (5,0)*+{_2}*\frm{o}="C";
 \endxy}\Ea
$$
Hence the differential $\delta_0$ acts only on input  white vertices of graphs by splitting each such white vertex
$
 \Ba{c}\resizebox{4mm}{!}{
 \xy
 (0,1.5)*+{_v}*\frm{o};
 \endxy
 }\Ea
$
into two new white vertices
$
 \Ba{c}\resizebox{10mm}{!}{ \xy
 (-5,2)*+{_{v'}}*\frm{o}="B";
  (5,2)*+{_{v''}}*\frm{o}="C";
 \endxy}\Ea
$
and redistributing all edges (if any) attached to $v$  in all possible ways among the new vertices $v'$ and $v''$. The cohomology
$$ \cE_1
\Def\left(\Ass_\infty \stackrel{\cF}{\lon} \f( \wh{\HoQP}_{0,1}^{\atop \circlearrowright})\right)[1]= H\left(\cE_0
\Def\left(\Ass_\infty \stackrel{\cF}{\lon} \f( \wh{\HoQP}_{0,1}^{\atop \circlearrowright})\right)[1], \delta_0\right)$$
is spanned by graphs all of whose white vertices are precisely univalent and skew symmetrized (see, e.g., Theorem 3.2.4 in \cite{Me-p} where this result is obtained from the cell complexes of permutahedra, or Appendix A in \cite{Wi1} for another purely algebraic argument) and hence is isomorphic (after erasing these no more needed  white vertices)  to $\Der(\wh{\HoQP}_{0,1}^{\atop \circlearrowright})$ as a graded vector space.
The boundary  condition (\ref{5: Boundary cond for formality map}) says that the induced differential $\delta_1$ in the complex
$
 \cE_1
\Def\left(\Ass_\infty \stackrel{\cF}{\lon} \f^{\geq 1}( \wh{\HoQP}_{0,1}^{\atop \circlearrowright})\right)[1]$
agrees precisely with the induced differential $d_1$ in
$\cE_1 \Der(\wh{\HoQP}_{0,1}^{\atop \circlearrowright})=\Der(\wh{\HoQP}_{0,1}^{\atop \circlearrowright})$ so that the morphism of the next pages of the spectral sequences,
$$
c_{\cF}^1: (\cE_1\Der(\wh{\HoQP}_{0,1}^{\atop \circlearrowright}), d_1)\lon
\left(\cE_1\Def(\cA ss_\infty \stackrel{\cF}{\rar} \f^{\geq 1}(\wh{\HoQP}_{0,1}^{\atop \circlearrowright}), \delta_1\right)
$$
is simply an isomorphism. By the comparison theorem, the morphism $c_{\cF}$ is a quasi-isomorphism.
\end{proof}

\subsection{Proof of Theorem {\ref{1: Class theorem of homog formalities}}}  Let $\cF_0$ be an arbitrary homogeneous formality map (in particular, the one constructed by M.Kontsevich in \cite{Ko2}). The above theorem implies
$$
H^{i+1}\left(\Def\left(\cA ss_\infty \stackrel{\cF_{0}}{\rar} \f^{\geq 1}(\wh{\HoQP}_{0,1}^{\atop \circlearrowright})\right)\right)=H^i(\fGC_2^{\geq 2}), \ \ \ \forall\ i\in \Z.
$$
In particular,
$$
H^1\left(\Def\left(\cA ss_\infty \stackrel{\cF}{\rar} \f(\wh{\HoQP}_{0,1}^{\atop \circlearrowright})\right)\right)=H^0(\fGC_2^{\geq 2}) =\grt
$$
To complete the  proof of Theorem {\ref{1: Class theorem of homog formalities}}  it remains to show that every infinitesimal deformation $\eta \in H^1\left(\Def\left(\cA ss_\infty \stackrel{\cF_0}{\rar} \f(\wh{\HoLB}_{0,1}^{\atop \circlearrowright})\right)\right)$ of any given homogeneous formality map $\cF_0$ exponentiates to a genuine
homogeneous formality map $\cF^\eta_0$, and then apply Lemma 3 from Thomas Willwacher paper \cite{Wi3}. The argument is standard and is based on the remarkable fact that the two dg Lie algebras
$$
\Der(\wh{\HoQP}_{0,1}^{\atop \circlearrowright})\ \text{and}\
 \Def(\wh{\HoQP}_{c,d}^{\atop\circlearrowright}\stackrel{\Id}{\rar}
\wh{\HoQP}_{c,d}^{\atop \circlearrowright})
$$
are identical -- after a degree shift --- as complexes
$$
\Der(\wh{\HoQP}_{0,1}^{\atop \circlearrowright})=
 \Def(\wh{\HoQP}_{c,d}^{\atop \circlearrowright}\stackrel{}{\rar}
\wh{\HoQP}_{c,d}^{\atop \circlearrowright})[1]
$$
but have really different Lie algebra structures (even, of different degrees). The element
$\eta$ (more precisely its arbitrary lift to a cycle) has degree zero when understood as an element of the dg Lie algebra $\Der(\wh{\HoQP}_{0,1}^{\atop \circlearrowright})$.
The Lie brackets of the latter respect the gradation by the total number of internal edges and legs
(which is always positive)
 so that it makes sense to consider its exponent $exp(\eta)$ as an element of the group $Aut(\wh{\HoQP}_{0,1}^{\atop \circlearrowright})$. Hence $exp(\eta)$  gives us a non-trivial Maurer-Cartan element of
the second dg Lie algebra
$$
\Def(\wh{\HoQP}_{c,d}^{\atop \circlearrowright}\stackrel{}{\lon}
\wh{\HoQP}_{c,d}^{\atop \circlearrowright})
$$
which is mapped into a non-trivial Maurer-Cartan element of the dg Lie algebra
$$
\Def\left(\cA ss_\infty \stackrel{\cF}{\rar} \f(\wh{\HoQP}_{0,1}^{\atop \circlearrowright})\right)
$$
giving us the required exponentiation of the infinitesimal deformation $\eta$ of $\cF_0$ to a genuine
formality morphism
$$
\cF_0^\eta: \cA ss_\infty {\lon} \f(\wh{\HoQP}_{0,1}^{\atop \circlearrowright})
$$
satisfying the boundary condition (\ref{5: Boundary cond for formality map}) and which differs at the infinitesimal level from $\cF_0$ precisely by the cycle $\eta$. The proof of part (i) is completed.

\sip

It is known \cite{Do,AM} that universal quantizations of generic Poisson structures are also classified by the same set of Drinfeld's associators. Hence if any two such quantizations agree on generic quadratic Poisson structures, then by the result in (i), they must correspond to the same associator and, therefore, must be equivalent. The proof of the theorem is completed.

\bip

\bip

\def\cprime{$'$}

\end{document}